\documentclass[12pt,a4paper]{article}
\usepackage{amssymb,amsmath,amsthm}
\usepackage[T1]{fontenc}
\usepackage{fullpage,parskip}
\usepackage{enumitem}
\usepackage{tabularray}
\usepackage[x11names]{xcolor}
\usepackage[table]{xcolor}
\usepackage{dsfont}

\usepackage{thmtools} % load correct theorem numbering stuff for cleverref
\usepackage[linktoc=all,hidelinks,colorlinks,unicode=true]{hyperref}
\usepackage[capitalise,compress,nameinlink,noabbrev]{cleveref} % Must be loaded after hyperref
\hypersetup{linkcolor={blue!70!black},citecolor={blue!70!black},urlcolor={blue!70!black}}

\declaretheorem[name=Theorem,numberwithin=section]{theorem}
\declaretheorem[name=Lemma,sibling=theorem]{lemma}

\declaretheorem[name=Corollary,sibling=theorem]{corollary}
\declaretheorem[name=Definition,sibling=theorem,style=definition]{definition}
\declaretheorem[name=Remark,sibling=theorem,style=definition]{remark}


\numberwithin{theorem}{section}

\DeclareMathOperator{\lcm}{lcm}
\DeclareMathOperator{\Id}{Id}
\DeclareMathOperator{\ord}{ord}
\DeclareMathOperator{\Res}{Res}
\DeclareMathOperator{\Var}{Var}
\DeclareMathOperator{\Bin}{Bin}

\newcommand\F{\mathbb{F}}
\newcommand\N{\mathbb{N}}
\newcommand\Z{\mathbb{Z}}
\newcommand\R{\mathbb{R}}
\newcommand\indic{\mathds{1}}
\newcommand\Prob{\mathbb{P}}
\newcommand\cA{\mathcal{A}}
\newcommand\cB{\mathcal{B}}

\newcommand\cH{\mathcal{H}}
\newcommand\cI{\mathcal{I}}

\newcommand\cO{\mathcal{O}}
\newcommand\cS{\mathcal{S}}
\newcommand\cT{\mathcal{T}}
\newcommand\cR{\mathcal{R}}
\newcommand\E{\mathbb{E}}
\newcommand\Q{\mathbb{Q}}
\newcommand\eps{\varepsilon}

\newcommand\bH{\mathbf{H}}

\title{Multiple coverings of the integers by arithmetic progressions}
\author{Paul Balister\footnote{Mathematical Institute, University of Oxford, United Kingdom,\\
\texttt{\{paul.balister, tommy.walkermackay\}@maths.ox.ac.uk}}\protect\footnotemark[1] \qquad
Tommy Walker Mackay\protect\footnotemark[1]\qquad}
\date{\today}

\begin{document}

\maketitle

\begin{abstract}
We say that a finite collection $\cA$ of arithmetic progressions is an $M$-cover for the integers if every integer is contained
in at least $M$ elements of~$\cA$. Let $M_\Z(s)$ be the largest $M$ such that an $M$-cover exists if each modulus
(i.e., common difference) of an arithmetic progression is allowed to occur in at most $s$ of the progressions in~$\cA$. We show that
$M_{\Z}(s)=e^\gamma s(\log s+(2 \log 2+o(1))\log\log s)$, where $\gamma$ is the Euler--Mascheroni constant.
We also prove generalisations of this result to rings of integers (and $S$-integers) in global fields, i.e., number fields or
global function fields, where we aim to cover the ring of integers by a finite collection of translates of ideals,
each ideal being used at most $s$ times. 
\end{abstract}

\section{Introduction} \label{sec:intro}

A \emph{covering system} of the integers is a finite collection of infinite arithmetic progressions (that is, subsets of $\Z$
of the form $a + m\Z$, where $m \in \N$ is the \emph{modulus} of the progression) whose union is~$\Z$.
Covering systems were first introduced by Erd\H os~\cite{erdos}, who used them to construct an infinite family
of positive odd numbers that cannot be written as the sum of a prime and a power of two, solving a problem of Romanov.
We say a covering system is \emph{distinct} or has \emph{multiplicity $1$}
if the moduli of the arithmetic progressions are distinct. Erd\H os asked whether for each $m > 0$ there exists a distinct
covering system in which each modulus is at least~$m$, a problem which was resolved in the negative by Hough~\cite{Hough} in 2015.
Other works (see e.g. \cite{BBMST21,Balister2022,CFT,houghnielsen}) have considered the existence of distinct covering systems with
additional constraints imposed on the moduli. In this paper, we relax the condition that the covering system is distinct
to allow each modulus to be used up to $s$ times, where $s$ is a fixed constant, and we will consider the circumstances
under which every integer appears in at least $M$ of the arithmetic progressions.

We extend the setting by considering covers of suitable rings $\cO$ by finite collections $\cA=\{a_i+I_i:i=1,\dots,n\}$
of cosets $a_i+I_i$ of ideals $I_i$ of~$\cO$. In the case when $\cO=\Z$,
these are just arithmetic progressions of modulus $m_i$ where $I_i=m_i\Z$. The most natural
extension is to the case when $\cO$ is a ring of integers (or more generally $S$-integers) in a global field.
However, it should be emphasised that we actually use very little number theory
in this paper, with all the key ideas being present in the case when $\cO=\Z$.

We start with a few definitions. We use the notation $[n]$ to denote $\{1, 2, \dots, n\}$ with the convention that $[0]=\emptyset$.

\begin{definition}[Ideal Multiplicity]
 Let $\cA = \{a_i+I_i:i\in[n]\}$ be a finite collection (multiset) of cosets of ideals of $\cO$.
 We define the \emph{ideal multiplicity}, or just \emph{multiplicity}, of $\cA$, $m(\cA)$, as the maximum number of times an ideal $I$ occurs
 as some $I_i$, i.e.,
 \[
  m(\cA) = \max_{I} \bigl|\{i \in [n]: I_i = I\}\bigr|.
 \]
\end{definition}

\begin{definition}[$M$-cover]
 Let $\cA = \{a_i+I_i:i\in[n]\}$ be a finite collection (multiset) of cosets of ideals.
 For $x \in \cO$, we say that $x$ is \emph{covered\/ \textup(exactly\textup) $M$ times by} $\cA$ if
 \[
  \bigl|\{i \in [n]: x \in a_i + I_i\}\bigr|= M.
 \]
 If all $x \in \cO$ are covered at least $M$ times by $\cA$, we call $\cA$ an \emph{$M$-cover} of $\cO$.
\end{definition}

\begin{definition}[$M_\cO(s)$]
 We define $M_{\cO}(s)$ to be the largest $M$ such that there exists an $M$-cover
 $\cA$ of multiplicity $m(\cA)\le s$. We set $M_{\cO}(s)=\infty$ if such $M$ are unbounded.
\end{definition}

In this paper we prove the following theorem.

\begin{theorem}\label{thm:main_theorem}
 Let\/ $\cO$ be the ring of\/ $S$-integers of some global field\/~$K$, where $S$ is a finite set of primes
 of the ring of integers of\/~$K$. Then
  \[
   M_{\cO}(s) = e^{\gamma} \alpha_{\cO} s \bigl(\log s + (2\log 2+o_\cO(1))\log \log s\bigr),
  \]
  where $\gamma \approx 0.5772$ is the Euler--Mascheroni constant and $\alpha_\cO$ is the residue
  of the pole at\/ $z=1$ of the corresponding zeta function\/ $\zeta_\cO(z)$.
\end{theorem}
We refer to the next section for details on the zeta function $\zeta_\cO(z)$ and the constant~$\alpha_\cO$,
but for now we just note that $\alpha_\cO=1$ when $\cO=\Z$. 

For the lower bound we use a fairly simple construction, using only ideals with prime factors of norm at most $s\log s$,
to obtain coverage up to $e^\gamma\alpha_{\cO}s(\log s+(1+o(1))\log\log s)$. This can be done by
packing cosets of ideals into cosets of ideals of smaller norm until we have
$e^\gamma\alpha_{\cO}s(\log s+(1+o(1))\log\log s)$ separate 1-covers. However, to obtain the correct factor
in the $\log\log s$ term, we use a randomised construction that provides extra coverage from cosets of ideals
containing one large prime factor of norm up to $s(\log s)^{2\log 2}$. This
randomised construction yields a cover that probably cannot be partitioned into 1-covers (see below),
but nevertheless uses these extra ideals very efficiently.

We note that any upper bound on $M_\cO(s)$ for $s=1$ and $\cO=\Z$ implies that there is an upper bound
on the minimum modulus of an Erd\H os covering system (or else for arbitrarily large $M$ an $M$-cover of
multiplicity $1$ could be constructed by taking a union of $M$ Erd\H os covering systems with distinct moduli,
each of which has minimum modulus greater than the maximum modulus of the previous covering system).
Making any progress on the upper bound on $M_{\Z}(s)$ is therefore necessarily at least as hard as showing that the
minimum modulus in any covering system is bounded, a result that was proved only in 2015 by Hough~\cite{Hough}.

For the upper bound on $M_{\cO}(s)$, we use the \emph{distortion method} of Balister, Bollob\'as, Morris, Sahasrabudhe and
Tiba~\cite{Balister2022}. We follow the description of this method as given in~\cite{BalisterBCC} with
a modification to take into account the need to show the non-existence of $M$-covers when $M>1$.
The idea is to show that ideals with prime factors of norm larger than $s(\log s)^{2\log 2}$ do not
contribute much mass when we weight the elements of $\cO/\lcm(I_1,\dots,I_n)$ appropriately, and so cannot
contribute much to the coverage. For ideals that have a prime factor with norm above $s(\log s)^2$, the distortion method
would apply in a manner similar to prior work. However, to get the  correct coefficient of $\log\log s$ in
\cref{thm:main_theorem}, we need to show none of the ideals with a prime factor of norm above $s(\log s)^{2\log 2}$
contribute significantly. For this, we combine the distortion method with a convexity argument to show
that one cannot get much coverage from ideals with just one prime factor of norm greater than $s(\log s)^{2\log 2}$. This convexity argument has been independently developed in a recent paper of Dai, Li, Li and Yi~\cite{dai2026moduli}.
The ideals with two or more large prime factors can then be handled using
a second moment argument, in a similar manner to previous applications of the distortion method.

As mentioned above, previous versions of the distortion method dealt with single covers.
One can easily modify these approaches to consider
$M$-covers that can be partitioned into $M$ single covers of $\cO$.
Indeed, we can adapt the method to consider coverings of $[M]\times\cO$ by sets of the form $\{n_i\}\times (a_i+I_i)$.
However, these `partitionable' $M$-covers are more restrictive than general $M$-covers.
For example, $\{\{1,2\}, \{1,3\}, \{2,3\}\}$ is a 2-cover of $\cS=\{1,2,3\}$, but no choice of $n_i$
gives a 1-cover of $[2]\times\cS$ by $\{\{n_1\}\times\{1,2\},\{n_2\}\times\{1,3\},\{n_3\}\times\{2,3\}\}$
(see also Example 1.1 of \cite{PS}).
If we define $M'_\cO(s)$ to be the corresponding maximum coverage by partitionable $M$-covers, then clearly
$M'_\cO(s)\le M_\cO(s)$, so our upper bound still applies to $M'_\cO(s)$. On the other hand, the
simple construction mentioned above (and described more precisely in \cref{sec:slb}) is partitionable.
Thus, we establish the following bounds on $M'_\cO(s)$:
\[
 e^{\gamma} \alpha_{\cO} s \bigl(\log s + (1+o(1))\log \log s\bigr)
 \le M'_{\cO}(s) \le e^{\gamma} \alpha_{\cO} s \bigl(\log s + (2\log 2+o(1))\log \log s\bigr).
\]
However, this is as far as we have been able to go with partitionable covers, and it
remains an open question to determine the correct coefficient in the $\log\log s$ term for $M'_{\cO}(s)$.

\section{Some number theoretic background}

We consider a global field $K$, i.e., a finite extension of either $\Q$ or some rational function
field $\F_q(t)$ over a finite field $\F_q$. Let $\cO_K$ be the usual ring of integers in $K$, i.e., the set
of elements of $K$ integral over $\Z$ or $\F_q[t]$ respectively. If $S$ is a finite set of primes of $\cO_K$,
set $\cO_{K,S}$ to be the set of $S$-integers in $K$, i.e., the set of elements $x\in K$ with $\ord_P(x)\ge0$ for
all $P\notin S$. We will consider covering systems in the case when $\cO$ is
of the form $\cO_{K,S}$ or $\cO_{K}$ (which is just $\cO_{K,\emptyset}$).
In the special case when $K=\Q$, we recover the original situation of the integers $\cO=\cO_{\Q}=\Z$.  

\begin{definition}[$\Id(\cO)$, $N(I)$] Write $\Id(\cO)$ for the set of non-zero ideals of $\cO$. For each $I\in\Id(\cO)$,
write $N(I):=|\cO/I|$ for the \emph{norm} of $I$, i.e., the number of residue classes mod~$I$.
\end{definition}

We note that $\cO$ is a \emph{residually finite Dedekind domain}:
\begin{itemize}
\item For every $I\in\Id(\cO)$, $N(I)$ is finite. Moreover, for any $x\in\R$
there are only finitely many ideals $I$ with $N(I)\le x$.
\item Every $I \in \Id(\cO)$ can be written uniquely (up to reordering) as $I = P_1^{a_1}P_2^{a_2}\dots P_k^{a_k}$,
where each $P_i$ is a prime ideal in $\cO$ and $a_i\in\N$.
\item The \emph{Chinese Remainder Theorem} holds: if 
$I = P_1^{a_1}P_2^{a_2}\dots P_k^{a_k}$ with $P_1,\dots,P_k$ distinct primes, then $\cO/I\cong \cO/P_1^{a_1}\times\cdots\times\cO/P_k^{a_k}$.
\item For a prime ideal $P$ and $v\ge1$, $N(P^v)=N(P)^v$. Indeed, if $\pi\in P\setminus P^2$ and $\cR\subseteq\cO$
is a set of representatives of $\cO/P$, then we can write any $x\in \cO/P^v$
uniquely in `base $\pi$' as $x=\sum_{i=0}^{v-1} a_i\pi^i+P^v$ with $a_i\in\cR$.
\end{itemize}

We remark that when considering $M$-covers of $\cO$, it is enough to consider cosets of \emph{non-zero}
ideals, since if $\cA$ is an $M$-cover then $\cA':=\{a+I \in \cA: I \ne 0\}$ is also an $M$-cover. Indeed, if
$\cA'$ were not an $M$-cover, then there would be some element $x\in\cO$ covered fewer than $M$ times by~$\cA'$.
But every element of the coset $x+Q$, $Q:=\lcm_{a+I\in \cA'}(I)=\bigcap_{a+I\in \cA'}I$,
is covered the same number of times by~$\cA'$,
so all of these elements are covered fewer than $M$ times by~$\cA'$. As $x+Q$ is an infinite set and each coset of the zero
ideal only covers a single element, this would contradict the fact that $\cA$ is an $M$-cover.

\begin{definition}[$\zeta_\cO$]
 We define the zeta function $\zeta_\cO(z)$ by
 \begin{equation}\label{e:zeta}
  \zeta_\cO(z) := \sum_{I\in\Id(\cO)}\frac{1}{N(I)^z}=\prod_{P\text{ prime}} \biggl(1-\frac{1}{N(P)^z}\biggr)^{-1}
 \end{equation}
 for $\mathrm{Re}(z) > 1$, and its analytic continuation elsewhere.
\end{definition}

It is known in the case when $\cO$ is the ring of integers of a global field that $\zeta_\cO$
is meromorphic on a disc centred at $1$ and that $\zeta_\cO(z)$ has a simple pole at $z=1$.
The same holds for $S$-integers: the primes of $\cO_{K,S}$ are in one to one correspondence with
the primes of $\cO_K$ that do not lie in $S$ (and with the same norms), so to obtain the zeta function for
$\cO_{K,S}$ one just removes the finite number of factors corresponding to $P\in S$ from the
Euler product in~\eqref{e:zeta}. As $\zeta_\cO(z)$ has a simple pole at $z=1$, we can make the following definition.

\begin{definition}[$\alpha_\cO$]
 \[
  \alpha_\cO := \Res_{z=1} \zeta_\cO(z).
 \]
\end{definition}

The constant $\alpha_\cO$ can be calculated from the analytic class number formula.
See \cite{neukirch}, Chapter VII \S 5 for number fields and \cite{Rosen1999} for function fields.
In particular, when $\cO = \Z$, $\zeta_\cO$ is the usual Riemann zeta function and $\alpha_{\cO} = 1$.
Another interesting special case is when $\cO = \F_q[t]$, in which case $\zeta_\cO(z)=\frac{1}{1-q^{1-z}}$
and $\alpha_{\cO} = \frac{1}{\log q}$.

\begin{definition}[$\pi_\cO(x)$] We define the prime counting function
\[
 \pi_\cO(x) := \bigl|\{P:P\text{ prime and }N(P)\le x\}\bigr|.
\]
\end{definition}

We shall use the following generalisation of (a weak form of) the prime number theorem. For number fields, see \cite{Landau1903},
and for function fields, our result is a direct consequence of Theorem 5.12 in~\cite{rosen} summed over all
sufficiently small degrees.

\begin{theorem}\label{thm:pnt_gen}
 \[
  \pi_{\cO}(x) = \Theta_{\cO}\Bigl(\frac{x}{\log x}\Bigr)\qquad\text{as}\quad x\to\infty.
 \]
\end{theorem}
Note that when $K$ is a number field, the stronger statement
\[
 \pi_{\cO}(x) \sim \frac{x}{\log x}
\]
holds, but when $K$ is a function field over $\F_q$, the $\Theta_{\cO}$ factor is needed to account
for the fact that prime ideals can only have norm equal to a power of~$q$.

\begin{definition}[$P_i$, $p_i$, $R_j(\cO)$]
As the ring $\cO$ has a countably infinite number of prime ideals, we can write these
ideals in non-decreasing order of norm as $P_1,P_2,\dots$. Write $p_i:=N(P_i)$, so that $p_1\le p_2\le\cdots$.
We define the finite products
\[
  R_j=R_j(\cO):=\prod_{i=1}^j \Bigl(1-\frac{1}{p_i}\Bigr)^{-1}.
\]
\end{definition}

We will use the following generalisation of Mertens' third theorem (a combination of Theorems 2 and 3 in \cite{Rosen1999}).
\begin{theorem}\label{thm:mertens_gen}
 \[
  R_j(\cO) = e^{\gamma} \alpha_{\cO} \log p_j  + O_{\cO}(1)\qquad\text{as}\quad j\to\infty.
 \]
\end{theorem}
Taking logarithms, it is easy to deduce the following strong, generalised form of Mertens' second theorem.
\begin{theorem}\label{thm:mertens_gen2}
 There exists a constant $C_{\cO}$ such that
 \[
  \sum_{i=1}^j\frac{1}{p_i}=\log\log p_j + C_{\cO}+O_{\cO}(1/\log p_j)\qquad\text{as}\quad j\to\infty.
 \]
\end{theorem}
%We also note that, as a simple consequence of \cref{thm:pnt_gen}, 
%\[
% \sum_{i=1}^\infty\frac{1}{p_i^2}=O_\cO(1).
%\]

\section{A simple lower bound}\label{sec:slb}

Suppose, as above, that the prime ideals of $\cO$ are listed in order of non-decreasing norm as $P_1,P_2,\dots$ and, for
all $j\in\N$ write $p_j:=N(P_j)$. Let $\cA=\{a_i+I_i:i\in[n]\}$ be a finite collection of cosets of ideals and
write $Q:=\lcm(I_1,\dots,I_n)=P_1^{\nu_1}\cdots P_J^{\nu_J}$. By the Chinese remainder theorem, we can identify
$\cS:=\cO/Q$ with $\cO/P_1^{\nu_1}\times\dots\times \cO/P_J^{\nu_J}$ and each $a_i+I_i\in\cA$ with a product
$A_i:=A_{i,1}\times\cdots\times A_{i,J}\subseteq \cO/P_1^{\nu_1}\times\dots\times \cO/P_J^{\nu_J}$, where
$A_{i,j}=(a_i+P_j^{v_{i,j}})/P_j^{\nu_j}$ and $I_i=P_1^{v_{i,1}}\cdots P_J^{v_{i,J}}$.
We note that $\cA$ is an $M$-cover if and only if the $A_i$, $i\in [n]$, form an $M$-\emph{cover} for~$\cS$,
i.e., every point $x\in\cS$ occurs in at least $M$ of the $A_i$. This geometric approach, considering
$M$-covers of $\cO$ as suitable coverings of product spaces with product sets, will be used in both the lower bound
and upper bound proofs of \cref{thm:main_theorem}. However, in the lower bound case we can further simplify the
picture by effectively removing all the higher powers of the~$P_j$.

So for this, we first introduce the notion of a \emph{discrete $M$-cover}.
Let $\cT := T_1 \times T_2 \times \dots \times T_J$
be a Cartesian product of finitely many finite non-empty sets. We say a set $H \subseteq \cT$ is a \emph{hyperplane}
if it is of the form $H := H_1 \times H_2 \times \dots \times H_J$, where for each $1 \leq j \leq J$, $H_j$
is either equal to $T_j$ or a singleton element of~$T_j$. We write $F(H) \subseteq [J]$ for the set of
\emph{fixed coordinates} of~$H$, i.e. the set $\{j \leq J: H_j \text{ is a singleton}\}$. (If $|T_j| = 1$ for any~$j$,
then we may choose arbitrarily whether or not $j \in F(H)$. In this case, a hyperplane is technically
a pair $(H,F(H))$, where $F(H)$ is some valid choice for the set of fixed coordinates of~$H$.)
We define a discrete $M$-cover and its multiplicity as follows.

\begin{definition}[Discrete $M$-cover of $\cT$]
 Let $\cT$ be a set as defined above, and let $\cH = \{H_i: i \in [n]\}$ be a finite collection (multiset)
 of hyperplanes in~$\cT$. For $x \in \cT$, we say $x$ is \emph{covered\/ \textup(exactly\textup) $M$ times by $\cH$} if
 \[
  |\{i \in [n]: x \in H_i\}| = M.
 \]
 If all $x \in \cT$ are covered at least $M$ times by $\cH$, we call $\cH$ an \emph{$M$-cover} of $\cT$.
\end{definition}

\begin{definition}[Multiplicity of a discrete $M$-cover]
 Let $\cH = \{H_i: i \in [n]\}$ be a discrete $M$-cover of $\cT$. We define the \emph{multiplicity} of $\cH$,
 $m(\cH)$, by
 \[
  m(\cH) = \max_{S \subseteq [J]} |\{i \in [n]: F(H_i) = S\}|.
 \]
\end{definition}

We shall write $M_\cT(s)$ for the largest $M$ such that there exists a discrete $M$-cover of $\cT$ with
multiplicity at most~$s$, and set $M_\cT(s) = \infty$ if such $M$ are unbounded.

We say that $\cT = T_1 \times T_2 \times \dots \times T_J$ is a \emph{discretisation}
of $\cO$ if for each $1 \leq j \leq J$, $|T_j| = p_j-1$. We first prove the following lemma relating the existence
of discrete $M$-covers to the existence of $M$-covers for~$\cO$.

\begin{lemma}\label{lem:discrete_cover}
 Let\/ $\cT$ be a discretisation of\/~$\cO$. Then $M_\cO(s) \geq M_{\cT}(s)$.
\end{lemma}
\begin{proof}
Suppose $\cH$ is an $M$-cover of $\cT$ with multiplicity~$s$. For each $j \leq J$,
pick some $\pi_j\in P_j\setminus P_j^2$ and identify $T_j$ with a set of representatives in $\cO$ of
the non-zero elements of $\cO/P_j$. This is possible as $|T_j|=p_j-1$. For each hyperplane
$H = H_1 \times \dots \times H_J \in \cH$ and for each function $g\colon F(H) \to \N$, we define
\[
 A(H,g) = \bigcap_{j \in F(H)}\bigl(h_j\pi_j^{g(j)-1}+P_j^{g(j)}\bigr),
\]
where $h_j \in T_j$ is the unique element such that $H_j = \{h_j\}$. Thus, by the Chinese remainder theorem,
$A(H,g)$ is a coset of the ideal $\prod_{j \in F(H)}P_j^{g(j)}$. Let $N$ be a large integer, to be chosen later,
and define $\cA$ to be the multiset given by
\[
 \cA := \{A(H,g): H \in \cH,\, g\colon F(H) \to [N]\}.
\]
Then $\cA$ is a finite collection of cosets of ideals that covers each element of $\cO \setminus \bigcup_{j=1}^{J} P_j^{N}$
at least $M$ times. Indeed, if $x\notin \bigcup_j P_j^N$ then, for each $j$, there exists some $g(j)\in[N]$ such
that $x\in P_j^{g(j)-1}$ but $x\notin P_j^{g(j)}$, so $x\in h_j\pi_j^{g(j)-1}+P_j^{g(j)}$ for some $h_j\in T_j$. 
But then $x$ lies in $A(H,g)$ for every $H$ that contains the point $(h_1,\dots,h_J)\in\cT$. As $\cH$ is an $M$-cover
of~$\cT$, $x$ is covered at least $M$ times by~$\cA$.

Furthermore, if $A(H,g)$ and $A(H',g')$ are cosets of the same ideal, then $g = g'$ and $F(H) = F(H')$.
Thus, as $\cH$ has multiplicity at most $s$, $\cA$ also has multiplicity at most~$s$.
    
We now construct another collection $\cA'$ of cosets of ideals which covers each element of
$\bigcup_{j=1}^{J} P_j^{N}$ at least $M$ times, has multiplicity $1\le s$, and which does not use
cosets of any of the same ideals as~$\cA$.
Pick any $K>J$ and recall that $p_K = N(P_K) = |\cO / P_K|$. Choose $N \geq MJp_K$, and let
$f\colon [M] \times [J] \times \cO/P_K \to [N]$ be an arbitrary injection. Define $\cA'$ by
\[
 \cA' = \bigl\{(a+P_K) \cap P_j^{f(m,j,a+P_K)} : 1\le m \leq M,\, 1\le j \leq J,\, a+P_K \in \cO/P_K\bigr\}.
\]
Then elements of $\cA'$ are cosets of distinct ideals $P_j^{f(m,j,a+P_K)}P_K$, all of which are distinct
from the ideals used by $\cA$ as none of those contain $P_K$ as a factor.
Also, for each $j$ and~$a$, $\cA'$ will cover the set $(a+P_K)\cap P_j^N$ at least $M$ times, and hence
$\cA'$ is an $M$-cover of $\bigcup_{j=1}^{J} P_j^{N}$.
Hence $\cA \cup \cA'$ is an $M$-cover of $\cO$ of multiplicity at most~$s$, demonstrating that $M_\cO(s) \geq M$.
As this holds for any $M$ such that there is an $M$-cover of $\cT$ of multiplicity at most~$s$, 
$M_\cO(s)\ge M_\cT(s)$.
\end{proof}

We now describe a simple method of constructing a set of hyperplanes with large coverage.

\begin{lemma}\label{l:split}
 Let $s \in \N$, let $T_1, T_2, \dots, T_J$ be finite non-empty sets and let $1 \leq j \leq J$. If\/ there exists a discrete $M$-cover of\/ $T_{j+1}\times\dots\times T_J$ of multiplicity
 at most $s+\lfloor s/|T_j|\rfloor$
 then there exists a discrete $M$-cover of\/ $T_j\times\dots\times T_J$ of multiplicity at most~$s$.
\end{lemma}
\begin{proof}
Suppose that we have such an $M$-cover $\cH$ of $T_{j+1}\times\dots\times T_J$. For each $S\subseteq\{j+1,\dots,J\}$
there are at most $s+\lfloor s/|T_j|\rfloor$ hyperplanes $H\in\cH$ with $F(H)=S$.
We now construct an $M$-cover $\cH'$ of $T_j\times\dots\times T_J$ of multiplicity at most~$s$.
For (up to) $s$ of the hyperplanes $H\in\cH$ with $F(H)=S$, include in $\cH'$ the hyperplane $H':=T_j\times H$ of
$T_j\times\dots\times T_J$ with $F(H')=S$. For the remaining at most $\lfloor s/|T_j|\rfloor$ hyperplanes
$H\in\cH$ with $F(H)=S$, include in $\cH'$ the $|T_j|$ distinct hyperplanes $H'_a:=\{a\}\times H$, $a\in T_j$,
which have $F(H'_a)=S\cup\{j\}$.
Note that, for a given $S\subseteq\{j+1,\dots,J\}$, there are now at most $s$ hyperplanes in $\cH'$ with
fixed coordinates $S$ and at most $s$ hyperplanes with fixed coordinates $S\cup\{j\}$.
Hence $\cH'$ has multiplicity at most~$s$.
As $\bigcup_{a\in T_j}H'_a=T_j\times H$, it is clear that $\cH'$ also
forms an $M$-cover of $T_j\times\dots\times T_J$.
\end{proof}

\begin{corollary}\label{c:lower_bound_gen}
 If\/ $\cT=T_1\times\cdots\times T_J$ is a discretisation of\/ $\cO$, then $M_\cT(s)\ge sR_J-O_\cO(J)$.
\end{corollary}
\begin{proof}
Let $s_0:=s$ and inductively define $s_j:=s_{j-1}+\lfloor s_{j-1}/|T_j| \rfloor$ for $j\ge1$.
For $j=J$ we have an $s_J$-cover of empty product $T_{j+1}\times\cdots\times T_J$ with multiplicity $s_J$
by simply taking $s_J$ copies of the trivial hyperplane with $F(H)=\emptyset$.
Then inductively for $j=J,J-1,\dots,1$, we obtain, by \cref{l:split}, an $s_J$-cover of $T_j\times\cdots\times T_J$
of multiplicity at most~$s_{j-1}$. Thus, we obtain an $s_J$-cover of $\cT$ of multiplicity at most~$s$. It remains to bound~$s_J$.

We have that $|T_j|=p_j-1$, so $s_j=s_{j-1}+\lfloor s_{j-1}/(p_j-1) \rfloor\ge (1-\tfrac{1}{p_j})^{-1}s_{j-1}-1$.
Thus, by induction,
\[
 s_J\ge s\prod_{j=1}^J\bigl(1-\tfrac{1}{p_j}\bigr)^{-1}-\sum_{j=1}^J\prod_{k=j+1}^J\bigl(1-\tfrac{1}{p_k}\bigr)^{-1}
 =sR_J-\sum_{j=1}^J\frac{R_J}{R_j}.
\]
Thus, it is enough to show that the last sum is $O_\cO(J)$. We split the sum into those with $p_j\le \sqrt{p_J}$ and those
with $p_j>\sqrt{p_J}$. For $p_j\le\sqrt{p_J}$ we just note that $R_j\ge 1$, so that by \cref{thm:pnt_gen,thm:mertens_gen},
\[
 \sum_{p_j\le\sqrt{p_J}}\frac{R_J}{R_j}\le \pi_{\cO}(\sqrt{p_J})\cdot R_J
 =O_{\cO}(\sqrt{p_J}/\log p_J\cdot \log p_J)=O_\cO(\sqrt{J}\log J)=o_\cO(J).
\]
For $p_j>\sqrt{p_J}$ we have
$R_J/R_j\le (e^\gamma\alpha_{\cO}\log p_J+O_{\cO}(1))/(e^\gamma\alpha_{\cO}\log\sqrt{p_J}+O_{\cO}(1))=O_{\cO}(1)$, and so
\[
 \sum_{\sqrt{p_J}\le p_j,\,j\le J}\frac{R_J}{R_j}=O_{\cO}(J),
\]
and the result follows.
\end{proof}

\begin{corollary}\label{c:lower_bound_simple}
 $M_{\cO}(s) \geq e^{\gamma} \alpha_{\cO} s (\log s + \log \log s) - O_{\cO}(s)$.
\end{corollary}
\begin{proof}
We take $J=s$, so $p_J=p_s=\Theta_\cO(s\log s)$, in \cref{c:lower_bound_gen} and apply \cref{lem:discrete_cover}
to obtain $M_\cO(s)\ge sR_s-O_\cO(s)$. The result then follows from \cref{thm:mertens_gen}. 
\end{proof}

\begin{remark} As noted in the introduction, this actually proves the existence of
such a cover that can be partitioned into $e^{\gamma} \alpha_{\cO} s (\log s + \log \log s) - O_{\cO}(s)$
$1$-covers of $\cO$.
\end{remark}

\section{An improved lower bound}

We note that if we take $J$ much larger than $s$ then for all $J > j$ much greater than $s$, $s_{j-1}<p_j-1$ in the proof of \cref{c:lower_bound_gen},
and so $s_{j-1}=s_j=s_{j+1}=\cdots=s_J$ and increasing $J$ no longer gains us any additional coverage via \cref{l:split}. Thus, there does not
seem to be a way of using ideals with prime factors of norm larger than about $p_s=\Theta_\cO(s\log s)$ by this method.
In this section we show that it is nevertheless possible to get slightly more coverage using cosets of
ideals with one large prime factor via a randomised approach.

We fix a small constant $\eps > 0$ and define
\begin{equation}\label{e:n}
 n:=\lfloor(\log s)^{2\log 2-2\eps}\rfloor.
\end{equation}
Let $J$ be maximal such that $p_J \leq ns$, let $J_1$ be maximal such that $p_{J_1} \leq s/(\log s)^{4\eps}$, and
let $J_2$ be maximal such that $p_{J_2} \leq s$. I.e.,
\begin{equation}\label{e:Js}
 J_1:=\pi_{\cO}(s(\log s)^{-4\eps}) < J_2:=\pi_\cO(s) <J:=\pi_\cO(ns).
\end{equation}
We also note that, as a consequence of \cref{thm:pnt_gen}, ratios of norms of consecutive primes are bounded above by a constant depending only on $\cO$, and so $p_{J_1}$, $p_{J_2}$ and $p_J$ are within an $O_\cO(1)$ factor
of $s/(\log s)^{4\eps}$, $s$, and $ns$ respectively.
Let $\cT:=T_1 \times T_2 \times \dots \times T_J$ be a discretisation of~$\cO$. We will construct an $M$-cover
of $\cT$ of multiplicity $s$, where $M = e^{\gamma} \alpha_O s(\log s + (2 \log 2 - o(1)) \log \log s)$.

By using hyperplanes $H$ with $F(H)\subseteq[J_2]$, we can construct, by \cref{c:lower_bound_gen}, an $M_1$-cover $\cH$
of multiplicity~$s$, where
\begin{equation}\label{e:M1}
 M_1=sR_{J_2}-O_\cO(J_2)=e^\gamma\alpha_{\cO}s\log s-O_{\cO}(s).
\end{equation}
To construct the rest of our cover, we use hyperplanes $H$ with fixed coordinate set $F(H) = S \cup \{j\}$,
where $S \subseteq [J_1]$ and $J_2 < j \leq J$. For each $J_2 < j \leq J$, we partition $T_j$
into $n$ sets of almost equal size:
$T_j = T_j^{(1)} \sqcup T_j^{(2)} \sqcup \dots \sqcup T_j^{(n)}$, where $|T_j^{(m)}|\le\lceil|T_j|/n\rceil$
for all $m\in[n]$. We will find it useful to group hyperplanes into collections that we will call \emph{cuboids}.

\begin{definition}[Cuboid]
 We say a collection $\bH = \{H_1, H_2, \dots, H_t\}$ of hyperplanes in $\cT = T_1 \times \dots \times T_J$
 is a \emph{cuboid} if there exists $S\subseteq[J_1]$, $J_2<j\le J$, $m\in[n]$ and choices of subsets $A_i\subseteq T_i$
 such that each $H_i$ is of the form
 \[
  H_i=A_1\times A_2\times\cdots\times A_{j-1}\times\{b_i\}\times A_{j+1}\times\cdots A_J,
 \]
 where $\{b_1,\dots,b_t\}=T^{(m)}_j$, $t=|T^{(m)}_j|$, and $F(H_i)=S\cup\{j\}$ (so $A_i=T_i$ when $i\notin S$ and
 $|A_i|=1$ when $i\in S$.)
\end{definition}

In other words, a cuboid is a collection of hyperplanes that are identical except in one fixed coordinate $j$,
where $J_2<j\le J$, and the values of this coordinate fill out one of the blocks $T_j^{(m)}$ of $T_j$.
Also, all other fixed coordinates must be at most~$J_1$.

Let $k_j:=\lfloor \frac{ns}{p_j}\rfloor-1$. Now for each $j$, $J_2<j\le J$, and for each $r\in [k_j]$, we construct 
a collection $\cA_r^j:=\{\bH_{j,r,S}:S\subseteq[J_1]\}$ of cuboids. To do this, choose a point
$U(j,r) = (U_1^{j,r},U_2^{j,r},\dots,U_{J_1}^{j,r})\in \cT_1:=T_1 \times T_2 \times \dots \times T_{J_1}$ uniformly at random,
independently for each pair $(j,r)$.
Now for each $S \subseteq [J_1]$, we pick $m=m_{j,r,S}\in [n]$ uniformly at random 
(independently of all other choices) and then take the cuboid $\bH_{j,r,S}$ that
covers $A_1 \times \dots \times A_{j-1} \times T_j^{(m)} \times A_{j+1} \times \dots \times A_J$,
where for each $i \in S$, $A_i = \{t_i\}$, and for each $i \notin S$, $A_i = T_i$.
Note that the multiplicity of the set of hyperplanes in $\bigcup_r\cA_r^j$ is at most
\[
 k_j\max_m|T^{(m)}_j|\le \bigl(\tfrac{ns}{p_j}-1\bigr)\bigl(\tfrac{p_j}{n}+1\bigr)
 = s-\tfrac{p_j^2-n^2s}{np_j}-1\le s,
\]
for sufficiently large $s$, where we have used that $p_j^2=s^{2+o(1)}>n^2s=s^{1+o(1)}$.
Write $\cH':=\bigcup_{j,r,S}\bH_{j,r,S}$ for the collection of hyperplanes in $\bigcup_{j,r}\cA_r^j$. 
As the set of fixed coordinates for hyperplanes in $\cA_r^j$ are distinct for different $j$,
$\cH'$ also has multiplicity at most~$s$. It remains to show that, with positive probability, $\cH'$ has high coverage.

We shall consider, for a fixed element $x:=(x_1,x_2,\dots,x_{J_1}) \in \cT_1$, the number of times that
$\{x_1\}\times\cdots\times\{x_{J_1}\}\times T_{J_1+1} \times \dots \times T_J$ is covered by $\cA_r^j$.
Let $Y = Y(j,r)$ be the number of $i \in [J_1]$ such that $U_i^{j,r} = x_i$. Then it can easily be seen that there are $Z = Z(j,r):=2^Y$
cuboids $\bH_{j,r,S}\in \cA_r^j$ such that the projection of $\bigcup \bH_{j,r,S}$ onto $\cT_1$ covers~$x$.
Indeed, these are just the ones where $S\subseteq\{i\in[J_1]:U_i^{j,r}=x_i\}$.
We show that $Z$ is concentrated in a particular sense that requires the following definition.

\begin{definition}[Weight-typical interval]
 Let $W$ be a random variable with finite expectation taking values in $\R_{\geq 0}$. We say an interval
 $I \subseteq \R_{\geq 0}$ is $\eps$-weight-typical for $W$ if $\E[W \cdot \indic(W \in I)] \geq (1-\eps) \E[W]$.
\end{definition}

\begin{lemma}\label{lem:Z_concentration}
 For sufficiently large $s$, $I:=[(\log s)^{2\log2-\eps}, (\log s)^{2\log2+\eps}]$
 is an $\eps$-weight-typical interval for~$Z$.
\end{lemma}
\begin{proof}
We have $Z=2^Y$ where $Y=\sum_{i=1}^{J_1}Y_i$ and
$Y_i:=\indic(U_i^{j,r}=x_i)$ are independent Bernoulli random
variables with $\Prob(Y_i=1)=\frac{1}{|T_i|}=\frac{1}{p_i-1}$.

Let $Y_1', Y_2', \dots, Y_{J_1}'$ be independent Bernoulli random variables with $\Prob(Y_i' = 1) = \frac{2}{p_i}$, and let $Y'=\sum_{i=1}^{J_1}Y_i'$.
Then, a direct calculation shows that for any $y_i\in\{0,1\}$,
\[
 \Prob(Y'_i=y_i)=\frac{\E[2^{Y_i}\indic(Y_i=y_i)]}{\E[2^{Y_i}]}.
\]
From this it follows, by independence, that for any choice of $y_1,\dots,y_{J_1}\in\{0,1\}$,
\[
 \Prob(Y'_1=y_1,\dots,Y'_{J_1}=y_{J_1})=\frac{\E[2^{Y}\indic(Y_1=y_1,\dots,Y_{J_1}=y_{J_1})]}{\E[2^Y]}.
\]
Hence, by linearity of expectation, it follows that for any $a<b$,
\[
 \Prob(Y'\in [a,b])=\frac{\E[2^{Y}\indic(Y\in [a,b])]}{\E[2^Y]}=\frac{\E[Z\cdot\indic(Z\in[2^a,2^b])]}{\E[Z]}.
\]
Thus, taking $[2^a, 2^b] = I := [(\log s)^{2\log2-\eps}, (\log s)^{2\log2+\eps}]$, it suffices to show that with probability at least $1-\eps$,
\[
 Y' \in [(2 - \eps/\log 2) \log \log s, (2 + \eps/\log 2) \log \log s].
\]
Now, by \cref{thm:mertens_gen2},
\[
 \E[Y'] = 2\sum_{i=1}^{J_1}\frac{1}{p_i}=2 \log \log (p_{J_1}) + O_\cO(1) = 2 \log \log s + O_\cO(1).
\]
Furthermore,
\[
 \Var(Y')=\sum_{i=1}^{J_1}\frac{2}{p_i}\biggl(1-\frac{2}{p_i}\biggr)\le \E[Y'].
\]
So, by Chebyshev's inequality, 
\[
 \Prob\bigl(Y' \notin[(1-\eps/2) \E(Y')), (1+\eps/2) \E(Y')]\bigr)\le \frac{\Var(Y')}{(\eps/2)^2\E[Y']^2}\to 0,
\]
as $s\to\infty$. This yields the desired result as, for sufficiently large~$s$,
\[
 [(1 - \eps/2) \E(Y')), (1+\eps/2) \E(Y')]\subseteq [(2 - \eps/\log 2) \log \log s, (2 + \eps/\log 2) \log \log s].\qedhere
\]
\end{proof}

We define $C_r^j(x)$ to be the number of times that $\{x_1\}\times\cdots\times\{x_{J_1}\}\times T_{J_1+1} \times \dots \times T_J$ is covered by~$\cA_r^{j}$, i.e., the minimum number of times any element of this set is covered.
We observe that there are $Z$ cuboids which could contribute to this cover, each of which covers
$\{x_1\}\times\cdots\times\{x_{J_1}\}\times T_{J_1+1} \times\cdots\times T_{j-1}\times T_j^{(m)}\times T_{j+1} \times\cdots \times T_J$
for exactly one
value of $m\in[n]$. Since this value of $m$ was selected uniformly at random for each of the $Z$ cuboids, the distribution
of $C_r^j(x)$ conditional on $Z = z$ is equal to the distribution of the minimum number of objects in any of $n$ boxes if $z$
objects are placed into one of $n$ boxes uniformly and independently at random. We now prove that if
$Z$ is significantly larger than~$n$, then with high probability $C_r^j(x)$ is not too much less than $Z/n$.

For this we shall use the following Chernoff bound.

\begin{theorem}\label{t:Chernoff}
 Suppose $X=\sum_{i=1}^nX_i$ is a sum of independent bounded random variables with $X_i\in[0,b]$. Then
 \[
  \Prob\bigl(X\le(1-\eps)\E[X]\bigr)\le \exp\bigl(-\eps^2\E[X]/2b\bigr).
 \]
\end{theorem}
\begin{proof}
The proof is standard, but we include it here
for completeness.

We first note that, by convexity of $\exp$, 
\[
 \E[e^{-\eps (X_i/b)}]\le 1+\E[X_i/b](e^{-\eps}-1)
 \le e^{(-\eps+\eps^2/2)\E[X_i/b]},
\]
where we have used $e^{-\eps}\le 1-\eps+\eps^2/2$
and $1+x\le e^x$ for $x,\eps\ge0$. Hence, by independence,
$\E[e^{-\eps (X/b)}]\le e^{(-\eps+\eps^2/2)\mu}$,
where $\mu:=\E[X]/b$. Now by Markov's inequality
applied to $e^{-\eps (X/b)}$ we have
\[
 \Prob\bigl(X\le(1-\eps)\E[X]\bigr)
 =\Prob\bigl(e^{-\eps(X/b)}\ge e^{-\eps(1-\eps)\mu}\bigr)
 \le \frac{\E[e^{-\eps(X/b)}]}{e^{-\eps(1-\eps)\mu}}\le e^{-\eps^2\mu+\eps^2\mu/2}
 =e^{-\eps^2\mu/2}.\qedhere
\]
\end{proof}

\begin{lemma}\label{lem:equidistribution}
 For any $\eps>0$ there exists a constant $K=K(\eps)$ such that if $z\ge Kn\log n$ and $z$
 objects are placed independently uniformly at random into $n$ boxes, then with probability at least $1-\eps$,
 every box contains at least $(1-\eps)z/n$ objects.
\end{lemma}
\begin{proof}
We may assume without loss of generality that $\eps<1$ and $n\ge2$. Define $K(\eps):=2\log(2/\eps)/(\eps^2\log 2)$.
The number of objects in any given box has binomial distribution $A \sim \Bin(z,1/n)$ which is the distribution
of a sum of $z$ independent Bernoulli random variables. By \cref{t:Chernoff} with $b=1$,
\[
 \Prob(A < (1-\eps)z/n) \le \exp\biggl(-\frac{\eps^2 z/n}{2}\biggr)
 \leq \exp\biggl(-\frac{\eps^2K}{2}\log n\biggr)=(\eps/2)^{\log_2 n}\le \eps/n,
\]
so, by the union bound, the probability that there exists a box with
fewer than $(1-\eps)z/n$ objects is no more than $\eps$, as desired.
\end{proof}

We can use this to bound $\E[C_r^j(x)]$. However, as we shall later be applying \cref{t:Chernoff}, we will
need to place an upper bound on the part of $C_r^j(x)$ that we shall be using. To this end, define
\begin{equation}\label{e:b}
 b:=(\log s)^{2 \log 2 + \eps} / n = (1+o(1)) (\log s)^{3\eps}
\end{equation}
and
\[
 \tilde C_r^j(x):=C_r^j(x)\cdot \indic(C_r^j(x)\le b).
\]

\begin{lemma}\label{lem:C_mean}
 We have, for sufficiently large $s$,
 \[
  \E[\tilde C_r^j(x)] \geq (1-3\eps)R_{J_1}/n.
 \]
\end{lemma}
\begin{proof}
Let $I$ be the interval from \cref{lem:Z_concentration} and note that, by \eqref{e:n} and \eqref{e:b}, we have
$I\subseteq [n(\log s)^\eps,nb]$. Since clearly $C_r^j(x)\le Z/n$ always holds,
$Z\in I$ implies that $C^r_j(x)\le b$. Hence, as $C_r^j(x)\ge0$,
\[
 \E[\tilde C_r^j(x)]=\E[C_r^j(x)\cdot\indic(C_r^j(x)\le b)] \geq \E[C_r^j(x)\cdot\indic(Z \in I)].
\]
As $n(\log s)^{\eps}\ge Kn\log n=n\cdot O(\log\log s)$ for large~$s$, we have, by \cref{lem:equidistribution}, that 
$C_r^j(x) \geq (1-\eps)Z/n$ with probability at least $1-\eps$ conditional on $Z \in I$.
Hence,
\[
 \E[\tilde C_r^j(x)]
 \geq \frac{(1-\eps)^2}{n} \E[Z \cdot \indic(Z \in I)].
\]
By \cref{lem:Z_concentration}, $I$ is an $\eps$-weight-typical interval for $Z$, so
$\E[Z \cdot \indic(Z \in I)] \geq (1-\eps) \E[Z]$. Thus,
\[
 \E[\tilde C_r^j(x)] \geq \frac{(1-\eps)^3}{n} \E[Z] \ge (1-3\eps)\E[2^Y]/n.
\]
But $\E[2^Y] = \prod_{i=1}^{J_1}\E[2^{Y_i}]=\prod_{i=1}^{J_1} (1 + \frac{1}{p_i-1}) = R_{J_1}$, so
$\E[\tilde C_r^j(x)]\geq (1-3\eps)R_{J_1}/n$.
\end{proof}

\begin{remark} We note that $R_{J_1}/n$ is much less than one, and by far the most common value for $C_r^j(x)$ is zero.
\end{remark}

Now each element of $\{x_1\}\times\cdots\times\{x_{J_1}\} \times T_{J_1+1} \times \dots \times T_J$ is covered at least
\[
 \sum_{j=J_2+1}^{J} \sum_{r=1}^{k_j} C_r^j(x)
\]
times by~$\cH'$. Our next lemma estimates the number
$k := \sum_{j=J_2+1}^{J} k_j$ of terms in this sum.

\begin{lemma}\label{lem:k_sum}
 For sufficiently large $s$,
 \[
  k \geq (2 \log 2-3\eps) ns\log \log s/\log s.
 \]
\end{lemma}
\begin{proof}
For $J_2<j\le J$ we have $k_j=\lfloor \frac{ns}{p_j}\rfloor-1\ge \frac{ns}{p_j}-2$.
Summing over $j$ gives
\[
 k\ge ns\sum_{j=J_2+1}^J\frac{1}{p_j}-2J=ns(\log\log p_J-\log\log p_{J_2}+O_\cO(1/\log p_{J_2}))-2J
\]
by \cref{thm:mertens_gen2}. Now by choice of $J$ we have $p_J\le ns$ and hence $J=O_\cO(ns/\log p_J)$.
Also, $\log p_J\sim \log p_{J_2}\sim \log s$ and $\log p_J-\log p_{J_2}=\log n+O_\cO(1)$.
Hence, using the inequality $\log x-\log y\ge \frac{x-y}{x}$ for $x>y>0$,
\[
 k\ge ns\frac{\log p_J-\log p_{J_2}}{\log p_J}+O_{\cO}(ns/\log p_J)\ge ns\frac{\log n+O_{\cO}(1)}{\log p_J}.
\]
Hence, by definition \eqref{e:n} of $n$, we have, for sufficiently large~$s$,
\[
 k\ge ns \frac{(2\log 2-2\eps)\log\log s+O_{\cO}(1)}{(1+o_\cO(1))\log s}
 \ge (2 \log 2-3\eps)ns \log \log s/\log s.\qedhere
\]
\end{proof}

Define
\[
 C(x):=\sum_{j=J_2+1}^J\sum_{r=1}^{k_j}C_r^j(x)
 \qquad\text{and}\qquad
 \tilde C(x):=\sum_{j=J_2+1}^J\sum_{r=1}^{k_j}\tilde C_r^j(x).
\]
Thus, $C(x)$ is the coverage of $\{x_1\} \times\cdots\times\{x_{J_1}\}\times T_{J_1+1} \times \dots \times T_J$
by $\cH'$ and $\tilde C(x)\le C(x)$.
We now give a lower bound on the expectation of $\tilde C(x)$ and hence also of $C(x)$.

\begin{lemma}\label{lem:big_C_mean}
 For sufficiently large $s$,
 \[
  \E[\tilde C(x)] \geq (2 \log 2 - 8\eps) e^\gamma \alpha_\cO s \log \log s.
 \]
\end{lemma}
\begin{proof}
From \cref{thm:mertens_gen} we have $R_{J_1}=e^\gamma\alpha_\cO \log p_{J_1}+O_\cO(1)=e^\gamma\alpha_\cO \log s+O_\cO(1)$,
so by \cref{lem:k_sum,lem:C_mean} we have
\begin{align*}
 \E[\tilde C(x)]&\ge (2 \log 2-3\eps) ns\log \log s/\log s\cdot (1-3\eps)(e^\gamma\alpha_\cO \log s+O_\cO(1))/n\\
 &\ge (2\log2 - 8\eps)e^\gamma\alpha_\cO s\log\log s.\qedhere  
\end{align*}
\end{proof}

We now use a Chernoff bound to bound the probability that $\tilde C(x)$ diverges from its mean by more
than a small constant factor. We will conclude that this probability is less than
$|\cT_1|^{-1}$, establishing by the union bound that, with positive probability, all of $\cT$
is covered at least $(2 \log 2 - O(\eps)) e^\gamma \alpha_\cO s \log \log s$ times by $\cA'$,
as desired. First, we bound the size of~$\cT_1$.

\begin{lemma}\label{lem:small_T_size}
 We have
 \[
  |\cT_1|=|T_1 \times T_2 \times \dots \times T_{J_1}| \leq \exp(O_\cO(p_{J_1})).
 \]
\end{lemma}
\begin{proof}
For each $i \leq J_1$ we have $|T_i| \leq p_{J_1}$, and by \cref{thm:pnt_gen}, we have that
$J_1 = O_\cO(p_{J_1}/\log p_{J_1})$. Hence,
\[
 \log |T_1 \times \dots \times T_{J_1}| \leq (\log p_{J_1})O_\cO(p_{J_1}/\log p_{J_1}) = O_\cO(p_{J_1}),
\]
as desired.
\end{proof}

\begin{lemma}\label{lem:lower_bound}
 $M_{\cO}(s)\ge e^\gamma\alpha_\cO s(\log s+(2 \log 2 - o_{\cO}(1))\log\log s)$.
\end{lemma}
\begin{proof}
As $\tilde C_r^j(x)\le b$, \cref{t:Chernoff} implies that
\[
 \Prob\big(\tilde C(x)\le (1-\eps)\E[\tilde C(x)]\bigr)
 \le \exp\bigl(-\eps^2\E[\tilde C(x)]/2b\bigr).
\]
Write $M_2:= (2\log 2-10\eps) e^\gamma \alpha_\cO s \log \log s$. We deduce from \cref{lem:big_C_mean} that
$M_2\le (1-\eps)\E[\tilde C(x)]$ and, as $C(x)\ge \tilde C(x)$, that
\[
 \Prob\big(C(x)\le M_2\bigr)\le\exp\bigl(-\Omega_{\cO}(\eps^2 s\log\log s/b)\bigr).
\]
By \cref{lem:small_T_size}, the number of possible values of $x$ is only $\exp(O_\cO(p_{J_1}))$,
so, by \eqref{e:Js} and~\eqref{e:b}, the expected number of values of $x\in\cT_1$ for which $C(x) \leq M_2$ is
\[
 \exp\bigl(O_\cO(s(\log s)^{-4\eps}) - \Omega_\cO(\epsilon^2 s (\log \log s)(\log s)^{-3\eps} )\bigr) < 1
\]
for sufficiently large~$s$. Hence, with positive probability, $\cH'$ is an $M_2$-cover for~$\cT$.
This can be combined with the $M_1$-cover $\cH$ from \eqref{e:M1} without increasing the multiplicity, as the cover from \eqref{e:M1} uses only ideals with all prime factors at most $P_{J_2}$. 
The result is an $M$-cover of $\cT$ of multiplicity at most $s$,
where $M := M_1+M_2\ge e^\gamma \alpha_\cO s(\log s + (2 \log 2 - 11 \eps) \log \log s)$.
Hence, $M_\cO(s) \ge M_\cT(s)\ge  e^\gamma \alpha_\cO s(\log s + (2 \log 2 - 11\eps) \log \log s)$.
Letting $\eps\to0$ then gives $M_\cO(s)\ge e^\gamma \alpha_\cO s(\log s + (2 \log 2 - o_\cO(1)) \log \log s)$, as desired.
\end{proof}

\begin{remark} We could define $J_2$ so that $p_{J_2}$ is anywhere between $p_{J_1}$ and $s\log s$.
Increasing $J_2$ would increase $M_1$, but would decrease $M_2$ by an essentially equivalent amount. 
\end{remark}

\section{The distortion method}

In this section, we describe the distortion method that we will use to bound $M_{\cO}(s)$ from above.
This method was originally introduced in~\cite{Balister2022} to give, among other results, a simplified solution to the minimum modulus problem, originally solved by Hough~\cite{Hough}, and has since been
used to study covering systems in more generality, including analogous results in number
fields~\cite{liwangyi}, polynomial rings over finite fields~\cite{LiWangWangYi2024} and global function fields~\cite{liwangwangyi,Wang}.
The approach used here will follow more closely the description of the method in~\cite{BalisterBCC},
where we inductively construct a non-uniform random variable that samples a point in $\cO/\lcm(I_1,\dots,I_n)$.
We then bound the expected number of times this point is covered. We introduce two new ideas: firstly, we extend the distortion
method to deal directly with $M$-covers rather than just $1$-covers, and secondly, we introduce a convexity
argument to bound a certain probability more tightly than is possible by using the second moment.

Let $\cA := \{a_i + I_i:i\in [n]\}$ be a collection of cosets of non-zero ideals, and let $Q := \lcm(I_1,\dots, I_n)$.
Let the prime factorisation of $Q$ be
\[
 Q = P_1^{\nu_1}P_2^{\nu_2} \cdots P_J^{\nu_J}.
\]
As before, we list all the primes in order of increasing norm and so set $\nu_i=0$ when the $i$th prime does not divide~$Q$. Let $p_i:=N(P_i)$ and
$\cS:=\cO/Q$, which we interpret via the Chinese Remainder Theorem as a Cartesian product
$\cO/P_1^{\nu_1}\times \cdots\times \cO/P_J^{\nu_J}$.
For convenience, let $N_j$ be the size of the $j$th factor, so $N_j:=N(P_j^{\nu_j})=p_j^{\nu_j}$.
We shall abuse notation slightly by identifying each $a_i+I_i\in\cA$ with the corresponding subset
$(a_i+I_i)/Q =\{x+Q:x\in a_i+I_i\}\subseteq \cS$.

Let $\cB_j:=\{i:P_j\mid I_i$ but $P_k\nmid I_i$ for $k>j\}$ be the set of indices $i$ where $P_j$
is the largest prime factor of $I_i$ (or to be more precise, the last prime factor of $I_i$ in the
given order $P_1, \dots, P_J$). For completeness, define $\cB_0:=\{i:I_i=\cO\}$.
Identifying $a_i+I_i$ with the corresponding subset of $\cO/P_1^{\nu_1}\times\cdots\times\cO/P_J^{\nu_J}$,
we see that, for $i\in\cB_j$, the condition $x=(x_1,\dots,x_J)\in a_i+I_i$ depends on $x_j$
but not on $x_k$ for $k>j$. Let $B_j:=\bigcup_{i\in\cB_j}(a_i+I_i)\subseteq\cS$ be the union of
the corresponding sets $a_i+I_i$.

For $x\in\cS$ and $1\le j\le J$, define
\[
 L_j(x):=\{(x_1,\dots,x_{j-1},y,x_{j+1},\dots, x_J):y\in\cO/P_j^{\nu_j}\}
\]
to be the set of points that differ from $x$ only in the $j$th coordinate. Equivalently,
$L_j(x)=x+Q_j\subseteq\cS$, where $Q_j:=\prod_{i\ne j}P_i^{\nu_i}$.
Define
\[
 \alpha_j(x):=\frac{1}{N_j}|L_j(x)\cap B_j|
\]
and
\[
 \beta_j(x):=\frac{1}{N_j}\sum_{i\in \cB_j}|L_j(x)\cap (a_i+I_i)|
\]
to be the proportion of coverage and average level of coverage respectively
of $L_j(x)$ by $\{a_i+I_i:i\in\cB_j\}$. Fix constants $\delta_1, \dots, \delta_J \in [0,\frac12]$ to be chosen later.

We now construct a random variable $X=(X_1,\dots,X_J)\in \cS$, where $X_j\in\cO/P_j^{\nu_j}$ for $j\ge1$.
The $X_j$ will not all be independent. Instead, we will bias the distribution of each $X_j$, based on the
values of $X_i$, $i<j$, in such a way so as to make $X$ avoid the set $B_j$ as much as possible, but with
the proviso that $X_j$ is not too far from uniform, even conditioned on the previous~$X_i$. The parameters $\delta_j$
determine how far from uniform $X_j$ is allowed to be by insisting that $X_j$ takes no value more than
$(1-\delta_j)^{-1}$ more often than a uniform random variable. More formally, we require that
\begin{equation}\label{e:defp}
 \Prob(X_j=x_j\mid X_1=x_1,\dots,X_{j-1}=x_{j-1})=
 \begin{cases}
  \frac{1}{N_j}\max\left\{0,\frac{\alpha_j(x)-\delta_j}{\alpha_j(x)(1-\delta_j)}\right\},&\text{if }x\in B_j;\\
  \frac{1}{N_j}\min\left\{\frac{1}{1-\alpha_j(x)},\frac{1}{1-\delta_j}\right\},&\text{if }x\notin B_j,
 \end{cases}
\end{equation}
(interpreted in the obvious limiting sense when $\alpha_j\in\{0,1\}$.)
We make a few comments about this definition. First note that $\alpha_j(x)$ does not depend on the
coordinates $x_j,\dots,x_J$, so the expression for the conditional probability is well-defined.
Second, we note that in both expressions, the first term inside the min or max is taken iff $\alpha_j(x)\le\delta_j$.
In particular, it is easy to check that \eqref{e:defp} defines a probability distribution.
Indeed, fixing all $x_i$, $i<j$, the number of choices of $x_j$ yielding
the first case is $\alpha_j(x)N_j$ by definition of $\alpha_j(x)$, and one can easily check that
\[
 \alpha_j(x)\max\left\{0,\frac{\alpha_j(x)-\delta_j}{\alpha_j(x)(1-\delta_j)}\right\}
 +(1-\alpha_j(x))\min\left\{\frac{1}{1-\alpha_j(x)},\frac{1}{1-\delta_j}\right\}=1.
\]
Hence, the sum of the conditional probabilities over $x_j$ is $1$ and so \eqref{e:defp}
inductively defines the joint distribution of $(X_1,\dots,X_j)$ in terms of the
joint distribution of $(X_1,\dots,X_{j-1})$.

Finally, we note that for any $x\in\cS$,
\begin{equation}\label{e:distort}
 \Prob(X_j=x_j\mid X_1=x_1,\dots,X_{j-1}=x_{j-1})\le \frac{(1-\delta_j)^{-1}}{N_j}.
\end{equation}

We now bound the average coverage of the random point $X$ in terms of $\beta_j$.

\begin{lemma}\label{lem:coverbound}
For $1 \leq j \leq J$, we have
\[
 \sum_{i\in\cB_j}\Prob(X\in a_i+I_i)\le \E\biggl[\max\biggl\{0,\frac{\beta_j(X)-\delta_j}{1-\delta_j}\biggr\}\biggr].
\]
\end{lemma}
\begin{proof}
By \eqref{e:defp} we have
\begin{align*}
 \sum_{i\in\cB_j}\Prob(X\in a_i+I_i&\mid X_1=x_1,\dots,X_{j-1}=x_{j-1})=\\
 &=\sum_{i\in\cB_j}|L_j(x)\cap(a_i+I_i)|\cdot \frac{1}{N_j}
  \max\biggl\{0,\frac{\alpha_j(x)-\delta_j}{\alpha_j(x)(1-\delta_j)}\biggr\}\\
 &=\beta_j(x)\max\biggl\{0,\frac{\alpha_j(x)-\delta_j}{\alpha_j(x)(1-\delta_j)}\biggr\}
 \le\max\biggl\{0,\frac{\beta_j(x)-\delta_j}{1-\delta_j}\biggr\},
\end{align*}
where in the last line we have used that $\alpha_j(x)\le \beta_j(x)$.
The result now follows by the law of total probability.
\end{proof}

To bound the expectation of terms involving $\beta_j(X)$, we will need to bound the probability that $X$ lies
in a coset of an arbitrary ideal. For this we shall use the following.

\begin{lemma}\label{lem:idealprob}
 For any coset $a+I$ of an ideal $I\mid Q$ \textup(which we identify with its image in $\cS=\cO/Q$\textup),
 \[
  \Prob(X\in a+I)\le \frac{1}{N(I)}\prod_{P_k\mid I}\frac{1}{1-\delta_k}.
 \]
\end{lemma}
\begin{proof}
We prove this by induction on the number of $k$ such that $P_k\mid I$. Clearly if no $P_k\mid I$ then $I=\cO$
and $\Prob(X\in a+I)=1=1/N(I)$. Now suppose $I=P_j^r I'$ where $P_j\nmid I'$ and $j$ is the largest index such that $P_j\mid I$.
Then $X\in a+I$ if and only if both $X\in a+I'$ and $X_j\in a+P_j^r$. But conditioning on $X\in a+I'$ only imposes
conditions on $X_i$ for $i<j$, so by \eqref{e:distort},
\[
 \Prob(X_j\in a+P_j^r\mid X\in a+I')
 \le \frac{|(a+P_j^r)/P_j^{\nu_j}|}{N_j}\cdot\frac{1}{1-\delta_j}
 =\frac{1}{(1-\delta_j)p_j^r}=\frac{N(I')}{(1-\delta_j)N(I)}.
\]
Hence, by induction,
\begin{align*}
 \Prob(X\in a+I)&=\Prob(X_j\in a+P_j^r\mid X\in a+I')\Prob(X\in a+I')\\
 &\le\frac{N(I')}{(1-\delta_j)N(I)}\cdot\frac{1}{N(I')}\prod_{P_k\mid I'}\frac{1}{1-\delta_k}
 =\frac{1}{N(I)}\prod_{P_k\mid I}\frac{1}{1-\delta_k}.\qedhere
\end{align*}
\end{proof}

We wish to bound $\beta_j(x)$. To do this, recall that
\[
 \beta_j(x) = \frac{1}{N_j} \sum_{i \in \cB_j} |L_j(x) \cap (a_i + I_i)|,
\]
where
\[
 \cB_j = \{i:P_j\mid I_i\text{ but }P_k\nmid I_i\text{ for }k>j\}.
\]
Now fix some $\eps>0$ and let $J_1$ be the largest index such that
\begin{equation}\label{e:J1}
 p_{J_1}\le s(\log s)^{2\log 2+\eps}.
\end{equation} 
 We shall define $\delta_j$ by
\begin{equation}\label{e:delta}
 \delta_j := \begin{cases}
  0,&\text{if } j\le J_1;\\
  \frac{1}{2},&\text{if } j >J_1.
 \end{cases}
\end{equation}
As a result, $(X_1,\dots,X_{J_1})$ is uniform on $\cO/P_1^{\nu_1}\times\cdots\times\cO/P_{J_1}^{\nu_{J_1}}$.
For $i\in\cB_j$, write $I_i=P_j^{r_i}I'_i$ where $P_j\nmid I'_i$. For $j>J_1$, we partition $\cB_j$ into two sets by defining
\[
\cB_j^0 := \{i \in \cB_j: P_k\mid I'_i\text{ implies } k\le J_1\},
\]
and $\cB_j^1 := \cB_j \setminus \cB_j^0$.
Define, for $d \in \{0,1\}$,
\[
 \beta_j^d(x) = \frac{1}{N_j} \sum_{i \in \cB_j^d} |L_j(x) \cap (a_i + I_i)|,
\]
so that $\beta_j(x) = \beta_j^0(x) + \beta_j^1(x)$.

\begin{lemma}\label{lem:alt_alpha_expression}
 For $d \in \{0, 1\}$, \[
  \beta_j^d(x) = \sum_{i\in\cB_j^d}\frac{1}{p_j^{r_i}}\indic(x \in a_i+I'_i).
 \]
\end{lemma}
\begin{proof}
By definition, $\beta_j^d(x)$ is the sum over $i\in\cB_j^d$ of the proportion of $y$ values where
$(x_1,\dots,x_{j-1},y,x_{j+1},\dots,x_J)\in a_i+I_i$. But this is zero if $x\notin a_i+I'_i$, and
$1/p_j^{r_i}$ when $x\in a_i+I'_i$.
\end{proof}

We now describe our strategy. Clearly, if $\cA$ is an $M$-cover, then
\[
 M\le \sum_i\Prob(X\in a_i+I_i)=\sum_{j=0}^J\sum_{i\in\cB_j}\Prob(X\in a_i+I_i),
\]
as the right-hand side is just the expected number of cosets in $\cA$ that cover~$X$.
We split the sum into three terms. The first is the sum over $j\le J_1$. For these values of $j$
we have $\delta_j=0$ and the distribution of $(X_1,\dots,X_j)$ is uniform on $\cO/P_1^{\nu_1}\times\cdots\times\cO/P_j^{\nu_j}$.
In this case we simply have 
\[
 \sum_{j=0}^{J_1}\sum_{i\in\cB_j}\Prob(X\in a_i+I_i)=\sum_{j=0}^{J_1}\sum_{i\in\cB_j}\frac{1}{N(I_i)}.
\]
For the second and third parts, we note that for $j>J_1$, $\delta_j=\frac12$. Thus
\begin{align*}
 \E\biggl[\max\biggl\{0,\frac{\beta_j(X)-\delta_j}{1-\delta_j}\biggr\}\biggr]
 &=\E\bigl[\max\{0,2\beta_j(X)-1\}\bigr]\\
 &\le \E\bigl[\max\{0,2\beta^0_j(X)-\tfrac12\}\bigr]+\E\bigl[\max\{0,2\beta^1_j(X)-\tfrac12\}\bigr]\\
 &\le \E\bigl[\max\{0,2\beta^0_j(X)-\tfrac12\}\bigr]+2\E\bigl[\beta^1_j(X)^2\bigr],
\end{align*}
where we have used convexity of $\max\{0,\cdot\}$ in the second line, and the inequality
$2x^2\ge \max\{0,2x-\frac12\}$ (which follows from rearranging $(2x-1)^2\ge0$) in the third.

Now $\beta^0_j(X)$ depends only on $(X_1,\dots,X_{J_1})$, which is uniform. This allows us to
use a convexity argument to find a configuration that maximises $\E[\max\{0,2\beta^0_j(X)-\tfrac12\}]$.
By analysing this particular configuration, we shall obtain a good upper bound on this term.
The last term, $2\E[\beta^1_j(X)^2]$, will be bounded using an explicit calculation
relying on \cref{lem:idealprob}. However, as all the terms in its expansion use at least two
large primes, it turns out that this term can also be shown to be very small.

To summarise, if we have an $M$-cover of $\cO$ of multiplicity $s$, then
\begin{equation}\label{e:split}
 M\le M_1+M_2+M_3,
\end{equation}
where
\begin{align*}
 M_1&:=\sum_{j=0}^{J_1}\sum_{i\in \cB_j}\frac{1}{N(I_i)},\\
 M_2&:=\sum_{j=J_1+1}^J \E\bigl[\max\{0,2\beta^0_j(X)-\tfrac12\}\bigr],\quad\text{and}\\
 M_3&:=\sum_{j=J_1+1}^J2\E\bigl[\beta^1_j(X)^2\bigr].
\end{align*}
The main term will be $M_1$, with both $M_2$ and $M_3$ minor error terms.
We bound each of these in the following sections.

\section{The upper bound: bounding \texorpdfstring{$M_1$}{M1}}

\begin{lemma}\label{lem:boundM1}
 For sufficiently large $s$ we have
 \[
  M_1=\sum_{j=0}^{J_1} \sum_{i \in \cB_j} \frac{1}{N(I_i)} \leq sR_{J_1}
  \le e^\gamma\alpha_\cO s(\log s+(2\log 2+2\eps)\log\log s).
 \]
\end{lemma}
\begin{proof}
There are at most $s$ copies of each ideal $I\mid P_1^{\nu_1}\cdots P_{J_1}^{\nu_{J_1}}$ occurring
as some $I_i$, $i\in\cB_j$, $j\le J_1$. Thus 
\begin{align*}
 \sum_{j=0}^{J_1} \sum_{i \in \cB_j}\frac{1}{N(I_i)}
 &\le s\sum_{v_1=0}^{\nu_1}\cdots\sum_{v_{J_1}=0}^{\nu_{J_1}}\frac{1}{N(P_1^{v_1}\cdots P_{J_1}^{v_{J_1}})}\\
 &= s\prod_{j=1}^{J_1}\sum_{v=0}^{\nu_j}\frac{1}{p_j^v}
 \le s\prod_{j=1}^{J_1}\Bigl(1-\frac{1}{p_j}\Bigr)^{-1}=sR_{J_1}.
\end{align*}
Now by \cref{thm:mertens_gen}, $sR_{J_1}=e^\gamma\alpha_\cO s \log p_{J_1}+O_\cO(s)$.
By \eqref{e:J1}, $p_{J_1}\le s(\log s)^{2\log 2+\eps}$, and so
$\log p_{J_1}\le\log s+(2\log 2+\eps)\log\log s$. Hence, for sufficiently large~$s$,
$sR_{J_1}\le e^\gamma\alpha_\cO s (\log s+(2\log 2+2\eps)\log\log s)$.
\end{proof}

\section{The upper bound: bounding \texorpdfstring{$M_2$}{M2}}

We start with a lemma that identifies a configuration that maximises the expression
$\E[\max\{0,2\beta^0_j(X)-\frac12\}]$.

Define a \emph{weighted collection of cosets}
$\cA=\{(a_i+I_i,w_i):i\in[n]\}$ to be a finite multiset
of cosets of ideals in $\cO$ together with non-negative
weights $w_i\in\R$. Define the \emph{coverage} $C_\cA(x)$ of $x\in\cO/Q$, where $Q:=\lcm(I_1,\dots,I_n)$, as
\[
  C_\cA(x):=\sum_{i\colon x\in a_i+I_i}w_i,
\] 
i.e., the sum of the weights of the cosets covering~$x$.

\begin{lemma}\label{lem:convexity_argument}
 Let\/ $\cA=\{(a_i+I_i,w_i):i\in[n]\}$ be a weighted collection of cosets and let
 $\cA_0:=\{(I_i,w_i):i\in [n]\}$ be the corresponding
 weighted collection of ideals.
 If\/ $f\colon\R\to\R$ is a convex function
 and $X$ is uniformly distributed on $\cO/\lcm(I_1,\dots,I_n)$, then
 \[
  \E\bigl[f(C_{\cA}(X))\bigr]\le \E\bigl[f(C_{\cA_0}(X))\bigr].
 \]
\end{lemma}
\begin{proof}
Write $Q:=\lcm(I_1,\dots,I_n)=P_1^{\nu_1}\cdots P_J^{\nu_J}$.
For each prime $P_j$, pick a set of representatives $\cR_j\subseteq\cO$ of $\cO/P_j$
and $\pi_j\in P_j\setminus P_j^2$. We can write any element
$a\in \cO/P_j^{\nu_j}$ uniquely in `base $\pi_j$' as a sum $a=\sum_{i=0}^{\nu_j-1}\alpha_i\pi_j^i+P_j^{\nu_j}$
where $\alpha_i\in \cR_j$. This allows us to construct a bijection
$\cO/P_j^{\nu_j}\to \cR_j^{\nu_j}=\cR_j\times\cdots\times \cR_j$ by mapping $a$ to $(\alpha_0,\dots,\alpha_{\nu_j-1})$.
A coset $a+P_j^v$, $v\le\nu_j$, corresponds to a hyperplane in $\cR_j^{\nu_j}$
in which the fixed coordinates form an initial segment of the coordinates of $\cR_j^{\nu_j}$.
Repeating this process for each prime, we can represent $\cO/Q$ as a product \[
\cT=\underbrace{\cR_1\times \cdots \times \cR_1}_{\nu_1}\times\cdots\times \underbrace{\cR_J\times \cdots \times \cR_J}_{\nu_J}.
\]
A coset $a_i+I_i$ then corresponds to a hyperplane
$H_i\subseteq \cT$. We can therefore translate the problem to the geometric case of discrete covers by hyperplanes. We will write $\cT = T_1 \times T_2 \times \dots \times T_K$, where $K = \sum_{j=1}^J \nu_j$ and $T_k = \cR_j$ for $\sum_{j'=1}^{j-1} \nu_{j'} < k \leq \sum_{j'=1}^{j} \nu_{j'}$.

So let $\cH=\{(H_i,w_i):i\in[n]\}$ be a weighted collection of hyperplanes of $\cT:=T_1\times\cdots\times T_K$,
where $H_i=H_{i,1}\times\cdots\times H_{i,K}$ with $H_{i,k}=\{h_{i,k}\}$ when $k\in F(H_i)\subseteq[K]$
and $H_{i,k}=T_k$ otherwise.
Fix a coordinate $k_0\in[K]$ and let $0\in T_{k_0}$ be some fixed chosen element of~$T_{k_0}$.
Define a weighted collection $\cH':=\{(H'_i,w_i):i\in[n]\}$ of hyperplanes, with
$H'_i=H'_{i,1}\times\cdots\times H'_{i,K}$, where $H'_{i,k}=H_{i,k}$ when either $k\ne k_0$ or $k=k_0\notin F(H_i)$, but
$H'_{i,k_0}=\{0\}$ when $k_0\in F(H_i)$.
For $x\in\cT$, define $C_{\cH}(x)=\sum_{i:x\in H_i}w_i$
and similarly for $C_{\cH'}(x)$. Let $X=(X_1,\dots,X_K)$ be a uniformly chosen random element of~$\cT$.
We shall show that $\E[f(C_{\cH}(X))]\le \E[f(C_{\cH'}(X))]$.

Fix $x = (x_1, \dots, x_K) \in \cT$, and let $I:=\{i\in[n]:x_k\in H_{i,k}\text{ for all }k\ne k_0\}$. Then if $X_k=x_k$ for $k\ne k_0$,
$C_{\cH}(X)=\sum_{i\in I:x_{k_0}\in H_{i,k_0}}w_i$ and similarly $C_{\cH'}(X)=\sum_{i\in I:x_{k_0}\in H'_{i,k_0}}w_i$.
Let $I':=\{i\in I:k_0\in F(H_i)\}$. Then
\[
 C_{\cH}(X)=\sum_{i\in I\setminus I'}w_i+\sum_{i\in I'}w_i\indic(X_{k_0}=h_{i,k_0}) =: g(X_{k_0})
\]
and
\[
 C_{\cH'}(X)=\sum_{i\in I\setminus I'}w_i+\sum_{i\in I'}w_i\indic(X_{k_0}=0) =: h(X_{k_0}).
\]
As $X_{k_0}$ is uniformly distributed in $T_{k_0}$,
\[
\E[g(X_{k_0})]=\E[h(X_{k_0})]=\sum_{i\in I\setminus I'}w_i+\sum_{i\in I'}w_i/|T_{k_0}|.
\]
However,
$h$ only takes the extreme values $\sum_{i\in I\setminus I'}w_i$ and $\sum_{i\in I}w_i$ whereas
$g$ always lies between these values. Hence, by Jensen's inequality
\begin{align*}
 \E[f(C_{\cH}(X))\mid X_k=x_k\text{ for }k\ne k_0]
 &=\E[f(g(X_{k_0}))]\\&\le \E[f(h(X_{k_0}))]
 =\E[f(C_{\cH'}(X))\mid X_k=x_k\text{ for }k\ne k_0].
\end{align*}
Thus, by the law of total probability, $\E[f(C_{\cH}(X))]\le \E[f(C_{\cH'}(X))]$.

Now we iterate this process. Define $\cH_0:=\cH$ and, for $k=1,\dots,K$,
let $\cH_k$ be the weighted collection of hyperplanes of $\cH_{k-1}$ but with the factors $H_{i,k}$
replaced by $\{0\}$ when $k\in F(H_i)$. Then, by the above, we have
$\E[f(C_{\cH}(X))]=\E[f(C_{\cH_0}(X))]\le \E[f(C_{\cH_1}(X))]\le\cdots\le \E[f(C_{\cH_K}(X))]$.
Now if $\cH:=\{(H_i,w_i):i\in [n]\}$ then $\cH_K=\{(H'_i,w_i):i\in [n]\}$
where $H_{i,k}=\{0\}$ when $k\in F(H_i)$ and $H'_{i,k}=H_{i,k}=T_k$ otherwise.
In the case of the weighted collection of cosets $\cA=\{(a_i+I_i,w_i):i\in [n]\}$, if
$\cH$ is the weighted collection of hyperplanes corresponding to $\cA$, then it is easy to see that
$\cH_K$ is the weighted collection of hyperplanes corresponding to $\cA_0=\{(I_i,w_i):i\in[n]\}$.
As the coverage functions $C_{\cA}(X)$ and $C_{\cA_0}(X)$ correspond directly to
$C_{\cH_0}(X)$ and $C_{\cH_K}(X)$ respectively, we have $\E[f(C_{\cA}(X))]\le \E[f(C_{\cA_0}(X))]$.
\end{proof}

\begin{corollary}\label{c:beta0}
 $\E[\max\{0,2\beta^0_j(X)-\frac12\}]$ is maximised
 when $a_i=0$ for all\/ $i\in\cB^0_j$.
\end{corollary}
\begin{proof}
Let $\cA:=\{(a_i+I'_i,1/p_j^{r_i}):i\in \cB^0_j\}$
Then by \cref{lem:alt_alpha_expression},
$\beta_j^0(x)=C_{\cA}(x)$. The result
follows from \cref{lem:convexity_argument} by
taking the convex function $f(x):=\max\{0,2x-\frac12\}$.
\end{proof}

We use this corollary to bound $\E[\max\{0,2\beta^0_j(X)-\frac12\}]$ above.

\begin{lemma} \label{lem:beta0_bound}
 For all real\/ $r\ge 1$,
 \[
  \E\bigl[\max\{0,2\beta_j^0(X) - \tfrac12\}\bigr]
  \le \frac12\Bigl(\frac{4s}{p_j-1}\Bigr)^r R_{J_1}^{2^r-1}.
 \]
\end{lemma}
\begin{proof}
By \cref{c:beta0} we may assume all $a_i=0$. Also, without loss of generality we may assume
that every ideal of the form $P_1^{v_1}\cdots P_{J_1}^{v_{J_1}}\cdot P_j^{v_j}$ with $0\le v_m\le \nu_m$ and $v_j\ge1$
occurs exactly $s$ times as some $I_i$, $i\in\cB_j^0$, as including more cosets can only increase $\beta_j^0$. Let $\cI_1$ be the set of ideals that are factors of 
$P_1^{\nu_1}\cdots P_{J_1}^{\nu_{J_1}}$ and let $C(x)$ be the number of ideals $I'\in \cI_1$
covering~$x=(x_1,\dots,x_J)\in\cS$.
Then, by \cref{lem:alt_alpha_expression},
\[
 \beta_j^0(x)=s\sum_{v=1}^{\nu_j}\sum_{I'\in \cI_1}\frac{1}{p_j^v}\indic(x\in I')
 =sC(x)\sum_{v=1}^{\nu_j}\frac{1}{p_j^v}\le \frac{sC(x)}{p_j-1}.
\]
Thus, for any $r\ge1$,
\begin{align}
 \E\bigl[\max\{0,2\beta_j^0(X) - \tfrac12\}]
 &\le \E\bigl[\max\{0,2sC(X)/(p_j-1) - \tfrac12\}]\notag\\
&\leq \frac{1}{2} \E\biggl[\frac{4sC(X)}{p_j-1} \indic\biggl(\frac{4sC(X)}{p_j-1} \geq 1\biggr)\biggr]\notag\\
&\leq \frac{1}{2} \E\biggl[\biggl(\frac{4sC(X)}{p_j-1} \biggr)^r\biggr]
 =\frac{1}{2}\Bigl(\frac{4s}{p_j-1}\Bigr)^r\E[C(X)^r].\label{e:rmoment}
\end{align}
It remains to bound $\E[C(X)^r]$. To do this, we note that
\[
 C(x)=\prod_{i=1}^{J_1}(1+v_i(x_i)),
\]
where $v_i(x_i)$ is the maximum $v\le \nu_i$ such that $x_i\in P_i^v$.
Indeed, $1+v_i(x_i)$ is the number of choices of $v_i$ such that $x\in P_i^{v_i}$.
Now for any non-negative integer~$v$,
\[
 (1+v)^r \le 2^r(2^r+1)\cdots(2^r+v-1)/v!.
\]
Indeed, this clearly holds for $v=0$, and for larger $v$ it holds by induction as
\[
 \frac{(1+v)^r}{(1+(v-1))^r}=(1+1/v)^r\le 1+(2^r-1)(1/v)=(2^r+v-1)/v,
\]
where the inequality is due to the convexity of $x^r$ for $r\ge1$ and the fact that $1/v\in[0,1]$.
Thus, as each $X_i$ is distributed uniformly in $\cO/P_i^{\nu_i}$ for $i\le J_1$,
\begin{align*}
 \E[(1+v_i(X_i))^r] &\le \sum_{v=0}^{\infty} \Bigl(1-\frac{1}{p_i}\Bigr)\frac{1}{p_i^v} \cdot (1+v)^r \\
 &\le \Bigl(1-\frac{1}{p_i}\Bigr)\sum_{v=0}^\infty \frac{1}{p_i^v}\cdot\frac{2^r(2^r+1)\cdots(2^r+v-1)}{v!}
 = \Bigl(1-\frac{1}{p_i}\Bigr)^{1-2^r}.
\end{align*}
Hence, as $X_1, X_2, \dots, X_{J_1}$ are independent,
\[
 \E[C(X)^r] \le \prod_{i=1}^{J_1}\Bigl(1-\frac{1}{p_i}\Bigr)^{1-2^r}= R_{J_1}^{2^r-1}.
\]
Substituting this into \eqref{e:rmoment} yields the desired result.
\end{proof}

We now sum these bounds for $J_1 < j \leq J$ to bound $M_2$.

\begin{lemma}\label{lem:boundM2}
 \[
  M_2=\sum_{j>J_1}\E\bigl[\max\{0,2\beta_j^0(X) - \tfrac12\}\bigr] = O_{\cO,\eps}(s).
 \]
\end{lemma}
\begin{proof}
By \cref{lem:beta0_bound}, we have that for any $r>1$,
\[
 \sum_{j>J_1}\E\bigl[\max\{0,2\beta_j^0(X) - \tfrac12\}\bigr]\le \frac{(4s)^r}{2} R_{J_1}^{2^r-1}\sum_{j>J_1}\frac{1}{(p_j-1)^r}.
\]
Now by \cref{thm:pnt_gen},
\[
 \sum_{j>J_1}\frac{1}{(p_j-1)^r}= O_{\cO,r}\bigl(p_{J_1}^{1-r}/\log p_{J_1}\bigr).
\]
Indeed, the sum over $j$ where $2^k p_{J_1}<p_j\le 2^{k+1} p_{J_1}$ is at most
\[
 2^{-kr}p_{J_1}^{-r}\cdot \pi_\cO(2^{k+1}p_{J_1})=O_{\cO}\bigl(2^{-kr}p_{J_1}^{-r}\cdot 2^{k+1}p_{J_1}/\log p_{J_1}\bigr)
 =2^{-(r-1)k}\cdot O_\cO(p_{J_1}^{1-r}/\log p_{J_1}),
\]
and summing over $k\ge0$ gives the result for $r>1$.
Fixing $r=1+\delta$, and writing $C_{\cO,r}$ for an appropriate constant depending only on $\cO$ and $r$,
\cref{thm:mertens_gen} implies that 
\begin{align*}
\sum_{j>J_1}\E\bigl[\max\{0,2\beta_j^0(X) - \tfrac12\}\bigr]
&\le C_{\cO,r} s^r (\log s)^{2^r-1} (s(\log s)^{2\log 2+\eps})^{1-r}/\log s\\
&\le C_{\cO,r}s (\log s)^{2^r-2+(2\log 2+\eps)(1-r)}\\
&\le C_{\cO,r}s(\log s)^{2^{1+\delta}-2-\delta(2\log 2+\eps)}.
\end{align*}
Now $2^{1+\delta}=2+(2\log 2)\delta+O(\delta^2)$,
so the exponent 
\[
 2^{1+\delta}-2-\delta(2\log 2+\eps)=O(\delta^2)-\eps\delta<0
\]
for sufficiently small $\delta$. In particular, we have
\[
 \sum_{j>J_1}\E\bigl[\max\{0,2\beta_j^0(X) - \tfrac12\}\bigr]=O_{\cO,\eps}(s)
\]
as required.
\end{proof}

\section{The upper bound: bounding \texorpdfstring{$M_3$}{M3}}

Recall that $\cB^1_j$ contains those $i\in\cB_j$ for which
the ideal $I_i$ contains at least one prime factor $P_{j'}\mid I_i$ with $J_1<j'<j$. As above, we write
$I_i=P_j^{r_i}I'_i$ where $P_j\nmid I'_i$.

\begin{lemma} \label{lem:Mj_bounds}
 With the above notation and $j>J_1$,
 \[ 
  \E[\beta_j^1(X)^2] \le \frac{s^2R_{J_1}^3}{(p_j-1)^2} \biggl(\Bigl(\frac{R_{j-1}}{R_{J_1}}\Bigr)^6-1\biggr).
 \]
\end{lemma}
\begin{proof}
By \cref{lem:alt_alpha_expression},
\[
 \E[\beta_j^1(X)^2]=\sum_{k,\ell\in\cB_j^1}\frac{1}{p_j^{r_k+r_\ell}}\Prob\bigl(X\in (a_k+I'_k)\cap(a_\ell+I'_\ell)\bigr).
\]
However, $(a_k+I'_k)\cap(a_\ell+I'_\ell)$ is either empty or is a coset of $I'_k\cap I'_\ell$, so by \cref{lem:idealprob},
\[
 \E[\beta_j^1(X)^2]\le\sum_{k,\ell\in\cB_j^1}\frac{1}{p_j^{r_k+r_\ell}}\cdot\frac{1}{N(I'_k\cap I'_\ell)}
 \prod_{P_m\mid I'_k\cap I'_\ell}\frac{1}{1-\delta_m}.
\]
Now $I_i$, $i\in\cB_j^1$, potentially runs over all ideals of the form $P_1^{v_1}\cdots P_{j-1}^{v_{j-1}}P_j^r$,
$0\le v_m\le \nu_m$, $1\le r\le \nu_j$, except those in which $v_m = 0$ for all $m>J_1$, with each such ideal
occurring $s$ times. If we let $\cI$ be the set of ideals that are factors of $P_1^{\nu_1}\cdots P_{j-1}^{\nu_{j-1}}$
and $\cI_1$ those that are factors of $P_1^{\nu_1}\cdots P_{j_1}^{\nu_{J_1}}$,
then
\begin{align*}
 \E[\beta_j^1(X)^2]
 &\le s^2\sum_{r,r'=1}^{\nu_j}\frac{1}{p_j^{r+r'}}
  \sum_{I,I'\in \cI\setminus\cI_1}\frac{1}{N(I\cap I')}\prod_{P_m\mid I\cap I'}\frac{1}{1-\delta_m}\\
 &\le s^2\biggl(\sum_{r=1}^{\nu_j}\frac{1}{p_j^r}\biggr)^2
 \sum_{I,I'\in \cI\setminus\cI_1}\frac{1}{N(I\cap I')}\prod_{P_m\mid I\cap I'}\frac{1}{1-\delta_m}\\
 &\le \frac{s^2}{(p_j-1)^2}\Biggl(\sum_{I,I'\in\cI}\frac{1}{N(I\cap I')}\!\prod_{P_m\mid I\cap I'}\frac{1}{1-\delta_m}-
  \sum_{I,I'\in\cI_1}\frac{1}{N(I\cap I')}\!\prod_{P_m\mid I\cap I'}\frac{1}{1-\delta_m}\Biggr).
\end{align*}
Note that we could also subtract terms corresponding to $I\in\cI_1$, $I'\in\cI\setminus \cI_1$,
and $I\in\cI\setminus \cI_1$, $I'\in\cI_1$, however it turns out that just subtracting terms with $I,I'\in\cI_1$
is enough to prove a sufficiently good bound.

We now consider the inner sums. As there are $2\nu+1$ values of $v,v'\in\{0,\dots,\nu_k\}$ with $\max\{v,v'\}=\nu\le \nu_k$,
we have
\begin{align*}
 \sum_{I,I'\in\cI}\frac{1}{N(I\cap I')}\prod_{P_m\mid I\cap I'}\frac{1}{1-\delta_m}
 &=\prod_{k<j}\Biggl(1+\frac{1}{1-\delta_k}\sum_{\substack{v,v'=0\\(v,v')\ne(0,0)}}^{\nu_k}\frac{1}{p_k^{\max\{v,v'\}}}\Biggr)\\
 &=\prod_{k<j}\Biggl(1+\frac{1}{1-\delta_k}\biggl(\frac{3}{p_k}+\frac{5}{p_k^2}+\cdots+\frac{2\nu_k+1}{p_k^{\nu_k}}\biggr)\Biggr)\\
 &\le\prod_{k<j}\biggl(1+\frac{1}{1-\delta_k}\frac{3p_k-1}{(p_k-1)^2}\biggr),
\end{align*}
and similarly for the sum over $I,I'\in\cI_1$. Hence, using that $\delta_k=0$ for $k\le J_1$ and $\delta_k=\frac12$ for $k>J_1$,
we obtain
\begin{align*}
 \E[\beta_j^1(X)^2]
 &\le\frac{s^2}{(p_j-1)^2}
 \biggl(\prod_{k<j} \Bigl(1+\frac{1}{1-\delta_k}\frac{3p_k-1}{(p_k-1)^2}\Bigr)-
 \prod_{k\le J_1} \Bigl(1+\frac{1}{1-\delta_k}\frac{3p_k-1}{(p_k-1)^2}\Bigr)\biggr)\\
 &=\frac{s^2}{(p_j-1)^2}
 \biggl(\prod_{k\le J_1} \Bigl(1+\frac{3p_k-1}{(p_k-1)^2}\Bigr)\biggr)
 \biggl(\prod_{J_1<k<j}\Bigl(1+\frac{2(3p_k-1)}{(p_k-1)^2}\Bigr)-1\biggr).
\end{align*}
Finally we note that
\[
 1+\frac{3p_k-1}{(p_k-1)^2} =\frac{1+1/p_k}{(1-1/p_k)^2}\le \Bigl(1-\frac{1}{p_k}\Bigr)^{-3},
\]
and
\[
 1+\frac{2(3p_k-1)}{(p_k-1)^2}\le \Bigl(1+\frac{3p_k-1}{(p_k-1)^2}\Bigr)^2\le \Bigl(1-\frac{1}{p_k}\Bigr)^{-6},
\]
so
\[
 \E[\beta_j^1(X)^2] \le
 \frac{s^2}{(p_j-1)^2}
 \cdot R_{J_1}^3\cdot\biggl(\Bigl(\frac{R_{j-1}}{R_{J_1}}\Bigr)^6-1\biggr),
\]
as required.
\end{proof}

\begin{lemma}\label{lem:boundM3}
 \[
  M_3=2\sum_{j>J_1} \E[\beta_j^1(X)^2] = o_{\cO}(s).
 \]
\end{lemma}
\begin{proof}
For any $0<a<b$ we have the inequality $\log b-\log a \le \frac{b-a}{a}$.
Thus, for $j>J_1$, \cref{thm:mertens_gen2} implies that
\[
 \log\frac{R_{j-1}}{R_{J_1}}=\log\log p_{j-1}-\log\log p_{J_1}+O_\cO(1/\log p_{J_1})
 \le\frac{\log p_{j-1}-\log p_{J_1}+O_{\cO}(1)}{\log p_{J_1}}.
\]
Hence, for $2^kp_{J_1}<p_j\le 2^{k+1}p_{J_1}$, 
\[
 \log\frac{R_{j-1}}{R_{J_1}}\le \frac{k+O_{\cO}(1)}{\log s}.
\]
Thus, if $k<\log s$, $\frac{R_{j-1}}{R_{J_1}}=1+O_{\cO}((k+1)/\log s)$ and so
\[
 \Bigl(\frac{R_{j-1}}{R_{J_1}}\Bigr)^6-1=O_{\cO}((k+1)/\log s).
\] 
On the other hand, if $k\ge \log s$ then \cref{thm:mertens_gen} implies that
\[
 R_{j-1}=e^\gamma\alpha_\cO(k\log 2+\log p_{J_1})+O_\cO(1)=O_\cO(k+\log s)=O_\cO(k),
\]
and
\[
  R_{J_1}=\Theta_{\cO}(\log p_{J_1})=\Theta_{\cO}(\log s),
\]
so 
\[
 \Bigl(\frac{R_{j-1}}{R_{J_1}}\Bigr)^6-1=O_{\cO}(k^6/(\log s)^6).
\]
Hence in general $(\frac{R_{j-1}}{R_{J_1}})^6-1=O_{\cO}((k^6+1)/\log s)$.
Now there are at most $\pi_{\cO}(2^{k+1}p_{J_1})\le O_{\cO}(2^{k+1}p_{J_1}/\log p_{J_1})$ primes
with $2^kp_{J_1}<p_j\le 2^{k+1}p_{J_1}$, so we have
\begin{align*}
 2\sum_{j>J_1} \E[\beta_j^1(X)^2]
 &\le 2\sum_{k=0}^\infty \frac{s^2R_{J_1}^3}{(2^kp_{J_1})^2}\cdot
  O_{\cO}\bigl((k^6+1)/\log s\bigr)\cdot O_{\cO}\bigl(2^{k+1}p_{J_1}/\log p_J\bigr)\\
 &\le \sum_{k=0}^\infty O_{\cO}(1)\cdot \frac{k^6+1}{2^k}\cdot \frac{s^2R_{J_1}^3}{p_{J_1}\log p_{J_1} \log s}\\
 &\le O_{\cO}(s^2 (\log p_{J_1})^2/p_{J_1} \log s)=o_{\cO}(s),
\end{align*}
where in the last line we have used that $p_{J_1}=\Theta_\cO(s(\log s)^{2\log 2+\eps})$.
\end{proof}

\section{Proof of \texorpdfstring{\cref{thm:main_theorem}}{Theorem \ref{thm:main_theorem}}}

\begin{lemma}\label{lem:upper_bound}
 $M_{\cO}(s)\le e^\gamma\alpha_\cO s(\log s+(2 \log 2 + o_{\cO}(1))\log\log s)$.
\end{lemma}
\begin{proof}
Suppose that we have an $M$-cover $\cA:=\{a_i+I_i:i\in [n]\}$ of multiplicity~$s$.
Then by \eqref{e:split} and \cref{lem:boundM1,lem:boundM2,lem:boundM3},
\[
 M\le M_1+M_2+M_3\le e^{\gamma}\alpha_{\cO}s(\log s+(2\log 2+2\eps)\log\log s)+O_{\cO,\eps}(s)+o_{\cO}(s).
\]
The result now follows on letting $\eps\to0$.
\end{proof}

Finally, \cref{thm:main_theorem} follows on combining the lower bound from \cref{lem:lower_bound}
with the upper bound from \cref{lem:upper_bound}.

\bibliographystyle{plain}
\bibliography{bibliography}

@incollection {BalisterBCC,
    AUTHOR = {Balister, Paul},
     TITLE = {{Erd\H os} covering systems},
 BOOKTITLE = {Surveys in combinatorics 2024},
    SERIES = {London Math. Soc. Lecture Note Ser.},
    VOLUME = {493},
     PAGES = {31--54},
 PUBLISHER = {Cambridge Univ. Press, Cambridge},
      YEAR = {2024},
      ISBN = {978-1-009-49053-5},
   MRCLASS = {05D10},
  MRNUMBER = {4749037},
}

@article {Balister2022,
    AUTHOR = {Balister, Paul and Bollob\'as, B\'ela and Morris, Robert and
              Sahasrabudhe, Julian and Tiba, Marius},
     TITLE = {On the {Erd\H os} covering problem: the density of the
              uncovered set},
   JOURNAL = {Invent. Math.},
  FJOURNAL = {Inventiones Mathematicae},
    VOLUME = {228},
      YEAR = {2022},
    NUMBER = {1},
     PAGES = {377--414},
      ISSN = {0020-9910,1432-1297},
   MRCLASS = {11B25 (11B05 11N36)},
  MRNUMBER = {4392459},
MRREVIEWER = {Donald\ Jason\ Gibson},
       DOI = {10.1007/s00222-021-01087-5},
       URL = {https://doi.org/10.1007/s00222-021-01087-5},
}

@article {BBMST21,
    AUTHOR = {Balister, Paul and Bollob\'as, B\'ela and Morris, Robert and
              Sahasrabudhe, Julian and Tiba, Marius},
     TITLE = {The {Erd\H os}-{Selfridge} problem with square-free moduli},
   JOURNAL = {Algebra Number Theory},
  FJOURNAL = {Algebra \& Number Theory},
    VOLUME = {15},
      YEAR = {2021},
    NUMBER = {3},
     PAGES = {609--626},
      ISSN = {1937-0652,1944-7833},
   MRCLASS = {11B25 (11A07 11N35)},
  MRNUMBER = {4261094},
MRREVIEWER = {Donald\ Jason\ Gibson},
       DOI = {10.2140/ant.2021.15.609},
       URL = {https://doi.org/10.2140/ant.2021.15.609},
}

@article {CFT,
    AUTHOR = {Cummings, Maria and Filaseta, Michael and Trifonov, Ognian},
     TITLE = {An upper bound for the minimum modulus in a covering system
              with squarefree moduli},
   JOURNAL = {Acta Math. Hungar.},
  FJOURNAL = {Acta Mathematica Hungarica},
    VOLUME = {175},
      YEAR = {2025},
    NUMBER = {1},
     PAGES = {1--25},
      ISSN = {0236-5294,1588-2632},
   MRCLASS = {11B25 (11A07)},
  MRNUMBER = {4880645},
MRREVIEWER = {Zhi-Wei\ Sun},
       DOI = {10.1007/s10474-024-01496-x},
       URL = {https://doi.org/10.1007/s10474-024-01496-x},
}

@article {erdos,
    AUTHOR = {{Erd\H os}, Paul},
     TITLE = {On integers of the form {$2^k+p$} and some related problems},
   JOURNAL = {Summa Brasil. Math.},
  FJOURNAL = {Summa Brasiliensis Mathematic\ae},
    VOLUME = {2},
      YEAR = {1950},
    NUMBER = {8},
     PAGES = {113--123},
      ISSN = {0039-498X},
   MRCLASS = {10.0X},
  MRNUMBER = {44558},
MRREVIEWER = {P.\ Scherk},
}

@article {Hough,
    AUTHOR = {Hough, Bob},
     TITLE = {Solution of the minimum modulus problem for covering systems},
   JOURNAL = {Ann. of Math. (2)},
  FJOURNAL = {Annals of Mathematics. Second Series},
    VOLUME = {181},
      YEAR = {2015},
    NUMBER = {1},
     PAGES = {361--382},
      ISSN = {0003-486X,1939-8980},
   MRCLASS = {11B25 (11A07 11N36)},
  MRNUMBER = {3272928},
MRREVIEWER = {Donald\ Jason\ Gibson},
       DOI = {10.4007/annals.2015.181.1.6},
       URL = {https://doi.org/10.4007/annals.2015.181.1.6},
}

@article {houghnielsen,
    AUTHOR = {Hough, Robert D. and Nielsen, Pace P.},
     TITLE = {Covering systems with restricted divisibility},
   JOURNAL = {Duke Math. J.},
  FJOURNAL = {Duke Mathematical Journal},
    VOLUME = {168},
      YEAR = {2019},
    NUMBER = {17},
     PAGES = {3261--3295},
      ISSN = {0012-7094,1547-7398},
   MRCLASS = {11B25 (11A07)},
  MRNUMBER = {4030365},
MRREVIEWER = {Donald\ Jason\ Gibson},
       DOI = {10.1215/00127094-2019-0058},
       URL = {https://doi.org/10.1215/00127094-2019-0058},
}

@article {Landau1903,
    AUTHOR = {Landau, Edmund},
     TITLE = {Neuer {B}eweis des {P}rimzahlsatzes und {B}eweis des
              {P}rimidealsatzes},
   JOURNAL = {Math. Ann.},
  FJOURNAL = {Mathematische Annalen},
    VOLUME = {56},
      YEAR = {1903},
    NUMBER = {4},
     PAGES = {645--670},
      ISSN = {0025-5831,1432-1807},
   MRCLASS = {99-04},
  MRNUMBER = {1511191},
       DOI = {10.1007/BF01444310},
       URL = {https://doi.org/10.1007/BF01444310},
}

@article {liwangwangyi,
    AUTHOR = {Li, Huixi and Wang, Biao and Wang, Chunlin and Yi, Shaoyun},
     TITLE = {On {Erd\H os} covering systems in global function fields},
   JOURNAL = {J. Number Theory},
  FJOURNAL = {Journal of Number Theory},
    VOLUME = {266},
      YEAR = {2025},
     PAGES = {269--280},
      ISSN = {0022-314X,1096-1658},
   MRCLASS = {11T55 (05B40 11A07 11Y35)},
  MRNUMBER = {4792479},
MRREVIEWER = {Donald\ Jason\ Gibson},
       DOI = {10.1016/j.jnt.2024.07.002},
       URL = {https://doi.org/10.1016/j.jnt.2024.07.002},
      NOTE = {arXiv:2402.03810},
}

@misc {liwangyi,
      title={On the minimum modulus problem in number fields}, 
      author={Huixi Li and Biao Wang and Shaoyun Yi},
      year={2023},
      eprint={2302.05946},
      archivePrefix={arXiv},
      primaryClass={math.NT},
      url={https://arxiv.org/abs/2302.05946}, 
      note={arXiv:2302.05946},
}

@book {neukirch,
    AUTHOR = {Neukirch, {J\"urgen}},
     TITLE = {Algebraic number theory},
    SERIES = {Grundlehren der mathematischen Wissenschaften},
    VOLUME = {322},
 PUBLISHER = {Springer-Verlag, Berlin},
      YEAR = {1999},
     PAGES = {xviii+571},
      ISBN = {3-540-65399-6},
   MRCLASS = {11Rxx (11-02 11S15 11S31 14C40)},
  MRNUMBER = {1697859},
MRREVIEWER = {Cornelius\ Greither},
       DOI = {10.1007/978-3-662-03983-0},
       URL = {https://doi.org/10.1007/978-3-662-03983-0},
}

@article {Rosen1999,
    AUTHOR = {Rosen, Michael},
     TITLE = {A generalization of {M}ertens' theorem},
   JOURNAL = {J. Ramanujan Math. Soc.},
  FJOURNAL = {Journal of the Ramanujan Mathematical Society},
    VOLUME = {14},
      YEAR = {1999},
    NUMBER = {1},
     PAGES = {1--19},
      ISSN = {0970-1249,2320-3110},
   MRCLASS = {11R47 (11R42 11R44)},
  MRNUMBER = {1700882},
MRREVIEWER = {A.\ Perelli},
}

@book {rosen,
    AUTHOR = {Rosen, Michael},
     TITLE = {Number theory in function fields},
    SERIES = {Graduate Texts in Mathematics},
    VOLUME = {210},
 PUBLISHER = {Springer-Verlag, New York},
      YEAR = {2002},
     PAGES = {xii+358},
      ISBN = {0-387-95335-3},
   MRCLASS = {11R58 (11R60 11T55)},
  MRNUMBER = {1876657},
MRREVIEWER = {Ernst-Ulrich\ Gekeler},
       DOI = {10.1007/978-1-4757-6046-0},
       URL = {https://doi.org/10.1007/978-1-4757-6046-0},
}

@article{LiWangWangYi2024,
  author  = {Huixi Li and Biao Wang and Chunlin Wang and Shaoyun Yi},
  title   = {On covering systems of polynomial rings over finite fields},
  journal = {Proceedings of the American Mathematical Society},
  volume  = {152},
  number  = {9},
  pages   = {3731--3742},
  year    = {2024},
  month   = sep,
  doi     = {10.1090/proc/16864},
  issn    = {0002-9939}
}

@article{Wang,
author = {Wang, Biao},
year = {2026},
number = {6},
pages = {1327--1336},
title = {Bounds for {Erdős} covering systems in global function fields},
volume = {22},
journal = {International Journal of Number Theory},
doi = {10.1142/S1793042126500703}
}

@article{PS,
author ={Hao Pan and Zhi-Wei Sun},
year = {2007},
volume = {135},
number = {11},
pages = {3515--3520},
title = {A Sharp Result on {$m$}-Covers},
journal = {Proceedings of the American Mathematical Society},
}

@article{dai2026moduli,
  title         = {On the Moduli of Distinct Covering Systems},
  author        = {Dai, Xiangjun and Li, Huixi and Li, Junfeng and Yi, Shaoyun},
  journal       = {arXiv preprint arXiv:2611.15110},
  year          = {2026},
  eprint        = {2611.15110},
  archivePrefix = {arXiv},
  primaryClass  = {math.NT},
  url           = {https://arxiv.org/abs/2611.15110}
}

\end{document}